\documentclass[journal]{IEEEtran}
\IEEEoverridecommandlockouts
\usepackage{cite}
\usepackage{amsmath,amsfonts,amssymb}
\usepackage{xcolor}
\usepackage{url}

\newcommand*{\QEDB}{\hfill\ensuremath{\square}}%

\ifCLASSINFOpdf
	\usepackage[pdftex]{graphicx}
 	\DeclareGraphicsExtensions{.pdf,.jpeg,.png,.eps}
\else
\fi
\usepackage{epstopdf}
\usepackage{ulem}
\usepackage{gensymb} % for the symbol "degree"
\usepackage{eso-pic}

\usepackage{tikz}
\usepackage{textcomp}
\usepackage{hyperref}
\usepackage{lipsum}

\newcommand\copyrighttext{% 
  \footnotesize \textcopyright 2017 IEEE. Personal use of this material is permitted.
  Permission from IEEE must be obtained for all other uses, in any current or future
  media, including reprinting/republishing this material for advertising or promotional
  purposes, creating new collective works, for resale or redistribution to servers or
  lists, or reuse of any copyrighted component of this work in other works.
  DOI: \href{http://ieeexplore.ieee.org/document/7815441/}{10.1109/TSP.2017.2652391}}
\newcommand\copyrightnotice{%
\begin{tikzpicture}[remember picture,overlay]
\node[anchor=south,yshift=1pt] at (current page.south) {\fbox{\parbox{\dimexpr\textwidth-\fboxsep-\fboxrule\relax}{\copyrighttext}}};
\end{tikzpicture}%
}

\begin{document}

\title{A study of  periodograms standardized using training data sets  and application to exoplanet detection}
\author{Sophia~Sulis, David~Mary and~Lionel~Bigot
\thanks{This work was supported by Thales Alenia Space, PACA region, CNRS project DETECTION/Imag'In and Programme National de Physique Stellaire (PNPS) of CNRS/INSU, France.}
\thanks{ The authors are with the  Universit\'e C\^ote d'Azur, OCA, CNRS, Laboratoire Lagrange, France
  (e-mail: Sulis.Sophia@oca.eu; David.Mary@unice.fr; Lionel.Bigot@oca.eu).}}

% ----------------------------------------- DOCUMENT  -----------------------------------------

\maketitle
\copyrightnotice

\begin{abstract}
	When the noise affecting time series is colored with unknown statistics, a difficulty for sinusoid detection  is to control the true significance level of the test outcome.
	This paper investigates the possibility of using training data sets of  the noise to improve this control. 
	Specifically, we analyze the performances of various detectors {applied to} periodograms standardized using training  data sets. Emphasis is put on sparse detection in the Fourier domain and on the limitation posed by  the necessarily finite size of the training  sets available in practice. We study  the resulting false alarm and detection rates and show that  standardization leads in some cases to powerful constant false alarm rate tests.
	The study is both analytical and numerical.   Although analytical results are derived in an asymptotic regime, numerical  results show that theory accurately describes the tests'  behaviour for moderately large sample sizes. 
	{Throughout the paper, an application of the considered periodogram standardization  is presented  for exoplanet detection in radial velocity data.}
\end{abstract}


\IEEEpeerreviewmaketitle

\section{Introduction}
\subsection{Considered problem and  application}
\label{secIA}
\IEEEPARstart{D}{etecting} sinusoids in noise is one of the most studied problems in signal processing. 
 Sinusoid detection is classically based on Fourier analysis, which benefits  from a 
 considerable set of statistical results. Often, the  assumptions required for these results to hold are :\\
{ (i)} Under the null hypothesis ($\mathcal{H}_0$ : noise only), 
the  random process $X$ is stationary Gaussian with known statistics;\\
{ (ii)} The time series $\{X(t_j)\}_{j=1,\cdots,N}$ is regulary sampled,  in which case Fourier analysis is performed through the 
 periodogram {\cite{Schuster_1898, Brillinger_1981}}
\begin{equation}
\label{eq1}
	P(\nu)   :=  \frac{1}{N\Delta t } \Big| \sum_{j=1}^{N} X(t_j)\mathrm{e}^{-{\rm{i}}2\pi\nu j\Delta t} \Big|  ^2,
\end{equation}
where $\nu$ is the considered frequency, $N$ the number of samples and $\Delta t$ the time sampling step;  \\
(iii) The number of samples $N$ is large (asymptotic regime);\\ 
(iv) The periodogram is evaluated at Fourier frequencies. \\
\indent Assumptions  (i)$-$(iv)  allow to characterize  the statistical properties of the periodogram as an  estimate of the power spectral density (PSD){{\cite{Bartlett_1950, Grenander_1957, Brillinger_1981, Priestley_1981, Brockwell_1991, Bloomfield_2000,Stoica_2005, Quinn_2001}}.
 In particular, any finite set of periodogram ordinates 
 near some fixed frequency 
 are asymptotically independent and exponentially distributed with parameters dependent on the noise PSD at the relevant frequencies (this assumption actually extends to linear processes other than Gaussian, see Th. 4 of \cite{Quinn_2001}).
In detection, such properties can be exploited
% often allows
%dealing with independent variables of known distributions allows  to  control 
to obtain the  probabilities of false alarm and detection of tests based on periodogram ordinates.}
% If  under the alternative hypothesis  ($\mathcal{H}_1$ : periodic signal in noise) the distribution of $X$ can also be precisely characterized, these properties translate into an accurate control of the probability of detection.

Despite a large literature on the subject, the detection of periodic signals remains an active field of research because in practical situations some (or all) assumptions (i)-(iv) above may not be met. Before reviewing some such situations, it is worth introducing the  application that motivated the present study %While this application, exoplanet detection,  will be used as an illustrative example throughout the paper,  
{ (the proposed detection approach may cover other applications, however).}

With almost { $3\,000$} confirmed exoplanets and nearly { $2\,500$} candidates (in { September} 2016,  \url{www.exoplanets.org}), the detection of extrasolar planets is an extremely active research field in Astrophysics since two decades. This field  benefits from constant technological improvements allowing extremely low noise detectors \cite{Pepe_2014b} and  long (months to years) spaceborne observations with high (a few tens of seconds) sampling rates{ \cite{Batalha_2014, Auvergne_2009, Rauer_2014}}. 
The Radial Velocity (RV)  technique is one method for exoplanet detection. When planets orbit a  host star, the resulting { gravitational} force creates a periodic displacement of the star. 
This induces a modulation of the RV of the star with respect to (w.r.t.) Earth, which translates into a periodic Doppler shift of the stellar light.  The RV  technique consists in detecting such variations in stellar RV time series {\cite{Fischer_2014,Perryman_2011}.}\\
%{ This induces a modulation of the stellar lines with respect to (w.r.t.) Earth, which is translated into a periodic Doppler shift of the stellar light and into a projected (radial) velocity. The RV  technique consists in detecting determinist periodic variations due to the presence of an orbiting object in the RV time series of the observed star {\cite{Fischer_2014,Perryman_2011}.}\\}
%
Modern instrumental performances have reached signal to noise ratios allowing, in principle, to detect exoplanets comparable to the Earth {by} the RV method%{ \cite{Pepe_2011}}
. However, at such low  levels of instrumental noise, a new critical and limiting issue appears. The stellar surface can be seen as a boiling fluid, with millions of convection cells generating  upward and downward plasma flows visible under the form of granules having typical lifetime of a few minutes. These motions generate random fluctuations in the measured RV of the star that can mimic, or hide, exoplanetary signatures.
%
%{ A current observational technique is to averaged several observations of a star in a night to reduce this convection effect (timescale of the granulation $\approx$ minutes/hours, supergranulation $\approx$ hours/days) but \cite{Meunier_2015} have shown that the resulting convective residual is still of some m/s which is the also the expected magnitude of Earth-like planets orbiting solar-like stars.}. 
%
The  case of  $\alpha$ Centauri Bb planet,  a detection claimed in 2012 \cite{Dumusque_2012} (with an evaluated $P$-value of $0.02 \%$) and since then subject to controversy \cite{Hatzes_2013, Rajpaul_2016}, is one example showing how our incomplete knowledge of the noise affects the reliability assigned to  exoplanet detection claims. This specific issue motivated the present study.\\
A key point is that in parallel to  technological advances,   astrophysicists have continuously improved stellar models and elaborated numerical simulation codes able to account for the complex interplay of various astrophysical processes in the star's interior and surface.  Recent works  demonstrate that 
{ granulation noise }  can be reproduced by large scale numerical simulations in a reliable way \cite{Bigot_2011}. This suggests the possibility of using such simulations to \textit{calibrate} the detection process, as RV time series are strongly affected by  this stellar noise. The present study shows that { one such calibration} indeed leads to improved  control of the statistical significance and to detection tests with increased power. 
%-----------------------------------------------------------------------------------------------------------------------------------------------------------------------------
\subsection{Unknown  noise statistics: related works}
\label{rel}
Turning back to the deviations encountered in practice w.r.t. assumptions (i)$-$(iv) above, this work will be primarily interested  in the case of
an incomplete knowledge of the noise statistics under the null hypothesis (\textit{i.e.}, relaxing condition { (i)}; we mention in the last section perspectives to relax { the} other conditions). 
In this situation, the distribution of $P$ under the null is not known and  consequently the significance level (the size) of the test is not known either. 
%Indeed, 
{  Constant false alarm rate (CFAR) detectors have been devised when the noise is white \cite{Fisher_1929,Chiu_1989,Shimshoni_1971}.}
%{ Of course},
%constant false alarm rate (CFAR) detectors,  \textit{i.e.}   tests robust to statistic mismatch under the null, exist (examples of such tests for sinusoid detection under  the assumption of an 
% independent and identically distributed (i.i.d.) Gaussian noise can be found in \cite{Fisher_1929,Chiu_1989,Shimshoni_1971}), but these tests  are  robust  against an unknown variance only if the noise is i.i.d. Gaussian. \\
When the noise is colored, the detection problem is  more complicated. In practice, two approaches can be followed. The first approach is simply to ignore possible noise correlations and to apply { tests designed for white noise}. However, as will be illustrated in this study, the statistical behavior of the resulting testing procedure may be hazardous, with unpredictable significance level and poor power (see 
%{ \cite{Sulis_2016a} and} 
also \cite{Kay_1999} on this point). \\
 A more sophisticated approach consists { of} estimating  the noise PSD (called $S_E$ below) so that condition { (i)} above is considered to hold approximately.  
This estimate  can then be used to calibrate the periodogram of the data $P(\nu)$, leading to a frequency-wise standardized periodogram of the form  
\begin{equation}
\widetilde{P}(\nu\;|\; \widehat{S}_E):=\frac{P(\nu)}{\widehat{S}_E(\nu)}.
\label{ptilde}
\end{equation}
Note that the classical Fisher's test \cite{Fisher_1929} {  standardizes the periodogram ordinates
% by the sample mean of the periodogram ordinates { 
%(that is, owing to Parseval identity in its discrete form, 
by the estimated  PSD of a white noise. Standardization \eqref{ptilde} can be seen as a generalization of this approach  (see \cite{Li_2014} for  a  recent review).} \\
  The estimate ${\widehat{S}_E(\nu)}$ can be  parametric 
  %(the PSD is considered to depend in a known manner on some parameters that are estimated) 
{  or non-parametric}.
  %(the PSD is non-parametric, but other parameters may be involved in the estimation procedure).  
%
Non parametric approches  originate from  seminal works  {  of  Whittle  \cite{Whittle_1952} and  Bartlett \cite{Bartlett_1955}}.
%,  where  ${\widehat{S}_E(\nu)}$ is  a truncated periodogram, and of  Bartlett \cite{Bartlett_1955}, who uses ``grouped periodogram''  assuming $S_E $ is approximately constant over some frequency bands. 
Parametric methods  often proceed by fitting  AutoRegressive (AR) or ARMA (AR Moving Average)  processes to the time series. 

{  A further complication arises  when multiple sinusoids are present under the alternative, as they perturb the estimation of the noise PSD
%. For instance, the performances of Fisher's test are well known to decrease when multiple sinusoids are present (in white noise) 
\cite{Siegel_1980,Priestley_1981}. Standardized tests  for this case  can be found in   \cite{Shimshoni_1971, Siegel_1980,Bolviken_1983a, Bolviken_1983b, Chiu_1989,Li_2014}.  
These tests are however non adaptive in the  number of sinusoids (which must be set {\textit{a priori}}) and designed for white noise. 
  For adaptive procedures for colored noise
see  Chap. 8 of \cite{Priestley_1981}, and 
\cite{Sachs_1993,Sachs_1994,Bhansali_1979,Truong_1990,Quinn_1991,Quinn_1999,Kavalieris_1994,Hannan_1961,Nicholls_1967}.
%regression techniques with robust M-estimators  \cite{Sachs_1993,Sachs_1994},   iterative procedures based on AR models of the periodogram residuals  once  the estimated harmonics are removed{\cite{Bhansali_1979,Priestley_1981,Truong_1990,Quinn_1991,Quinn_1999}},
 % model selection procedures based on  information theoretic criteria \cite{Kavalieris_1994}, attenuation of the influence of ordinates possibly contaminated  by sinusoids through modified window periodograms  \cite{Hannan_1961,Nicholls_1967} or double window periodograms  (Chap. 8 of \cite{Priestley_1981}). 
Techniques reducing the influence of signal peaks under the alternative are proposed in  \cite{Chiu_1990,Gryca_1998}. 
%peak-insensitive Yule-Walker spectral estimates are proposed in  \cite{Chiu_1990} and ARMA-like models are considered in \cite{Gryca_1998}. 
Different approaches,  related to standardization \eqref{ptilde}, can be found in  \cite{White_1999,Lu_2005,Liavas_1998,Zheng_2012}.
%use ratios of { AR} spectral estimates  \cite{White_1999}, {of pairs of periodogram ordinates for symmetry testing purposes   \cite{Lu_2005}, of periodogram ordinates over distinct frequency bands for detecting a signal with specified frequency {\cite{Liavas_1998},} or use} a geometric mean (instead of  arithmetic mean in Fisher's test)  for standardization \cite{Zheng_2012}.
\\ When following Generalized Likelihood Ratio (GLR) approaches { for detecting  multiple sinusoids in unknown number,  the GLR  must} be combined with model selection procedures. While  sharp model selection criteria and CFAR detectors exist under white noise, the correlated case remains an open problem \cite{Nadler_2011}.}
  %%%%%%%%%%%%%%%%%%
  In the particular field of exoplanet detection  {using RV}, we find similar families of techniques
  {  \cite{Koen_2015a,Koen_2015b,Hatzes_2013,Tuomi_2012}}.  
 % . In the white noise case, several works use standardized test 
  %in the wake of Fisher's approach \cite{Koen_2015a,Koen_2015b}. To account for unknown colored noise, 
  % tests for planet detection are in practice done after ``whitening'' the noise through  preprocessing stages (local trend filtering or 
   % harmonic removal), which exploit
   % \textit{a priori} knowledge about known or presumed noise sources ({\textit{e.g.}}  instrumental drifts or phenomena driven by stellar rotation%, { like spots}
%) \cite{Hatzes_2013}. Parametric noise  models are also used, as ARMA processes in \cite{Tuomi_2012} or more specific empirical  models of stellar PSD. 
%\cite{Harvey_1985,Dumusque_2011}. 
%{ (e.g; Harvey's laws \cite{Harvey_1985,Dumusque_2011}). }

In conclusion, regarding the problem of assessing tests' significance levels for multiple sinusoids detection in noise, an inspection of the literature shows that:\\
%\begin{itemize}
%	\item 
$\bullet$	For white noise of unknown variance, several studies provide accurate results for standardized test statistics of the form \eqref{ptilde}, \textit{e.g.},    \cite{Shimshoni_1971, Siegel_1980,Bolviken_1983a, Bolviken_1983b, Chiu_1989}.\\
	$\bullet$ For colored noise with unknown { PSD}, we are not aware of works studying the  false alarm rate when AR/ARMA or other models are used for test standardization as in \eqref{ptilde}. The difficulty in this case is the dependence of the distribution of $\widehat{S}_E$ on  estimated   {
	parameters, which complicates} the analytical characterization of the  distribution of $\widetilde{P}$.\\
The procedures described in \cite{Hannan_1961, Bhansali_1979,Priestley_1981}  provide asymptotic control of the false alarm rate. These procedures do not operate explicitly on test statistics of the form \eqref{ptilde} but rather on windowed periodograms that depend on several parameters. We found that {  these tests are  in practice sensitive to  parameter}  setting
{ and that
%. While it may be possible to set up automatic procedures for tuning these parameters,
estimating these parameters  impacts the significance level at which the tests are conducted}. 
This  level  can of course be approximated by simply neglecting  { the influence of 
such} a `preprocessing stage' { (dealing, \textit{e.g.}, with model order selection, filtering, adaptive window design, or standardization)}. For instance,  { we might pretend that $\widehat{S}_E={S}_E$ in \eqref{ptilde}.
% { at least asymptotically if $\widehat{S}_E$ is consistent estimate}. 
As  will be highlighted in  Sec. \ref{sec6} and \ref{sec7}, the actual significance level obtained when doing so may  however be} far from the assumed one. This leaves open the question of designing both powerful and CFAR tests for unknown colored noise and we propose  such tests in this paper. 
% \end{itemize}

Before closing this literature survey, we mention a few tests designed for a particularly interesting configuration of the detection problem,   
which is the so-called \textit {rare and weak} setting. In this setting, the sinusoids 
are both of small amplitudes { w.r.t. the noise level} and in small (and unknown) number w.r.t. the number of samples $N$. 
%This  is a difficult problem, which is particularly relevant  
%in the framework of Earth-like planet detection in RV data.
{   When viewed in the Fourier domain, sinusoid detection can } be casted as a sparse heterogeneous mixtures problem, which has attracted much attention in the last decade {\cite{Donoho_2004, Ingster_2010, Walther_2011, Moscovich_2016, Gontscharuk_2016}}. { We will see that while  fixed and adaptive  (in the number of sinusoids) procedures  lead to inconsistencies  when  noise correlations are ignored (because then the statistics under the null hypothesis are wrongly specified), such tests keep their nominal properties, with the CFAR property added, when applied to periodograms  suitably standardized with training data sets. In the particular case of the Higher Criticism \cite{Donoho_2004}, standardization \eqref{ptilde} is  an alternative approach
 to that of \cite{Hall_2010}.}   
%-----------------------------------------------------------------------------------------------------------------------------------------------------------------------------
\subsection{This study}
The present study  (an extended version of {\cite{Sulis_2016a}}) focuses  on the statistical characterization of test
statistics when both the PSD of the colored  noise and  the parameters of the sinusoids are unknown. We propose a detailed analysis of the effects of periodogram standardization by means   of, say, $L$ training time series{, which are}
% (with $L\ll N$), {\textit{i.e.}},
 independent realizations of the noise process alone,
and of the gain  that can be expected by using such training signals in a detection framework. \\
We, of course, make the important assumption that such a training data set is  available. Beyond the case of exoplanet detection  considered here,
one may imagine various situations where training signals can be obtained. In astronomical instruments for instances, secondary optical paths are often devoted to monitor `empty'  regions of the sky or calibration stars
{
%, in order to provide samples of the noise sources alone (a typical case is in the detection of
%stellar oscillations 
\cite{Gupta_2001}. \\}
Note that training noise vectors are routinely used for detection in radar systems, with however, an important difference.  
Adaptive test statistics   in radar typically use estimates of the covariance matrix of the training  vectors 
and therefore require $L>N$ for this matrix to be nonsingular. 
%\cite{Trees_2002}. 
This is a very different regime from that considered here, where $L\ll N$.\\
In the present study, we assume that the training data set is unbiased, in the sense that an averaged periodogram obtained from an infinitely large batch would converge uniformly to the true noise PSD. In practice, finite (possibly small) batch sizes can be encountered. For this reason, we say below that the noise is partially unknown and we address the effect of {  small values of $L$} on the detection performances. \\
Because one important objective of this study is to obtain analytical characterization of the test performances, we consider here a regular sampling. Comments on how to relax this assumption are discussed at the end of the paper.
Also,  our results are asymptotic in the number of samples $N$, which is characteristic of time series analysis. However, we will also pay attention to whether asymptotic theory accurately describes reality for finite sample sizes { through simulations}.\\
{ In the considered application framework of exoplanet detection  in RV data the working hypotheses are justified because accurate   simulations of stellar noise can be produced to form training data sets}. These  simulations are however computationally demanding. Obtaining a simulation of 100 days, for a star similar to the one shown  in \cite{Bigot_2011}, takes about 3 months of computing time on 120 cores
on modern clusters. 
Consequently, realistic values of $L$ are in the range of one to, say,  a hundred at most. This motivates the study of the impact of estimation noise in the proposed standardization approach. \\
We proceed as follows. { Sec.} \ref{sec2} presents the model and the detection approach. Sec. \ref{sec3} recalls classical results regarding periodogram's distribution.
Sec. \ref{sec4} and  \ref{sec5}  derive false alarm and detection rates for several tests. Sec. \ref{sec7}
is a numerical study. Table \ref{tab1} summarizes the main notations used in the paper.
  {
\begin{table}[ht!] 
 \caption{Table of notations}
   \label{tab1} 
\begin{tabular}{|c|c|} 
\hline
	$  N$ & Number of data points \\%[0.1cm] 
%	\hline
%	$t_j$ & Observation date (regular sampling) \\
%	[0.1cm] 
%	& such as $\Delta t = Cste$  \\[0.1cm] 
\hline
	$  X(t_j)$ &  Regularly sampled time series  \\[0.1cm] 
\hline
	$ E(t_j)$ &  Zero-mean stationary Gaussian colored noise  \\[0.1cm] 
\hline
	$\nu,  \nu_k$ &  Continuous frequency, Fourier frequency\\[0.1cm] 
	\hline
	$ S_E(\nu)$ &  Noise PSD  \\[0.1cm] 
\hline
	$ r_E$ &  Noise autocorrelation function \\[0.1cm] 
\hline
	 $  N_s, \alpha_q, f_q, \varphi_q $ &  Parameters of model \eqref{hyp}: Number of sines,\\
	 & sines' amplitude, frequency and phase \\[0.1cm] 
\hline
	$ N_p $ & Number of exoplanets orbiting target star \\[0.1cm] 
\hline
	$\textsf{T}_p, K_p,M_p$ & Planet period and its six Keplerian parameters\\
	$e_p,\omega_p,t_0,\gamma_0$ &  \\%[0.1cm] 
\hline
	$ N_C $ & Proxy for $N_s$ \\[0.1cm] 
\hline
	$ P(\nu)$ &  Classical periodogram \\[0.1cm] 
\hline
	 $  L$ &  Number of available training data sets  \\[0.1cm] 
\hline
	$ \overline{P}_L(\nu)$ &  Periodogram averaged with $L$ training data sets \\[0.1cm] 
\hline
	$ \widetilde{P}(\nu|\overline{P}_L)$ &  Periodogram standardized by $ \overline{P}_L$ \\[0.1cm] 
\hline
	$ \Omega $ & Indices set of considered Fourier frequencies \\[0.1cm] 
\hline
	$ \lambda_k := \lambda(\nu_k)$ & Non centrality parameter \\[0.1cm] 
\hline
	$ F_{\lambda_k}(d_1,d_2) $ & Non central Fisher-Snedecor distribution \\
	& with $d_1$ and $d_2$ degrees of freedom \\[0.1cm] 
%\hline
%	$  \varphi_{F}(\gamma,d_1;d_2)  $  & Density and CDF  \\
%	$ \Phi_{F}(\gamma,d_1,d_2)$& of a F distribution. \\[0.1cm
\hline
	$Z,z $ & Scalar random variable, one realization of $Z$ \\[0.1cm] 
\hline
	$ v_Z(z) $ & Pr$(Z>z)$ (observed p-value)  \\[0.1cm] 
\hline
	$ V_{Z} $ & P-value as a random variable ($ V_{Z}\sim {{\cal{U}}_{[0\; 1]}}) $ \\
%	& $v_Z$ is one realization of $V_{Z}$ \\[0.1cm] 
%	& in this case, are used for random quantities. \\[0.1cm] 
\hline
	${\bf{Z}}=[Z_1,\hdots,Z_N]^\top $ & Vector of random variables \\[0.1cm] 
\hline
	$ V_{{\bf{Z}},k} $ & P-value associated to component $Z_k$ of  ${\bf{Z}}$\\[0.1cm] 
\hline
	$ V_{{\bf{Z}},(k)} $ & k-th ordered p-value of  {\bf{Z}} \\[0.1cm] 
%\hline
%	$\rm T_M(Z) $ & Test of the maximum $Z$ value \\[0.1cm] 
%\hline
%	$\rm T_C(Z) $ & Test of the $\rm N_C^{\text{th}}$ maximum $Z$ value \\[0.1cm] 
%\hline	
%	$\rm T_F(Z) $ & Fisher test  \\[0.1cm] 
%\hline
%	$\rm HC(Z) $ & Higher Criticism test \\[0.1cm] 
%\hline
%	$\rm BJ(Z) $ & Berk-Jones test \\[0.1cm] 
%\hline
%	$\rm P_{FA} $ & False alarm probability \\[0.1cm] 
%\hline
%	$\rm P_{DET} $ & Detection probability \\[0.1cm] 
%\hline
%	$\rm I $ & Regularized incomplete Beta function \\[0.1cm] 
%\hline
%	$\rm \gamma $ & Detection threshold \\[0.1cm] 
%\hline
%	$\rm \widehat{S}_{E,AR} $ & Noise PSD estimated by AR process (estimated order  \\
%	& $\rm \hat{o}_{AR}$, coefficients $\rm c_j$  and predicted error variance $\rm \hat{\sigma}^2$). \\[0.1cm] 
\hline
\end{tabular}
\label{default}
\end{table}%
}
% \indent Addressing the question of the impact of estimation noise in standardized tests first requires,  for the two hypotheses of our model (Sec. \ref{sec2}), to  investigate  the  statistics of the classical, averaged  and standardized periodograms. This is the purpose of  Sec. \ref{sec3}. The second step (Sec. \ref{sec4}) is to select several tests  for which the benefits gained from the proposed standardization can be highlighted and quantified.  We opt for a sample of  tests covering the different cases of  single and multiple sinusoids detection.  At this stage we will be in position to analyze  the statistics of these tests when { applied to} periodograms standardized by means of a simple  non parametric method (namely, an averaged periodogram of the training data sets). Sec. \ref{sec5}  derives false alarm and detection rates for some detectors and highlights their CFAR property.  For comparison,  tests { applied to} periodograms standardized using  classical parametric methods will be considered in Sec. \ref{sec6}. The last part of the study is a numerical evaluation of the theoretical results and a comparison of the different methods  (Sec. \ref{sec7}). 
%Sec. \ref{conc} summarizes the main contributions and discusses some perspectives. 
 %
 %-----------------------------------------------------------------------------------------------------------------------------------------------------------------------------
 %-----------------------------------------------------------------------------------------------------------------------------------------------
\section{Statistical model and detection approach}
 \label{sec2}
%%%%%%
We consider the two hypotheses: 
\begin{equation} 
  \left\{         
      \begin{aligned}
	 \text{ ${\cal{H}}_0$ : } {X}(t_j) &= {\displaystyle{{{E}}}}(t_j) \\
	 \text{ ${\cal{H}}_1$ : } X(t_j) &=  \sum_{q = 1}^{N_s} \alpha_q \sin(2\pi f_q t_j+\varphi_q)+ E(t_j) \\
      \end{aligned}
    \right.
    \label{hyp}
\end{equation}
where $X(t_j=j\Delta t),\; j=1,\hdots,N$, is an evenly sampled data time series, $E(t_j)$ is a zero-mean second-order stationary Gaussian noise, 
{ with $\inf (S_E(\nu) ) > 0$ and ${ \sum_{u \in \mathbb{R}}} | r_E(u) | < \infty$}.   This is the noise of which we assume a set of training time series is available. { To simplify the presentation, we consider for the rest of the paper a unit  sampling step $\Delta t=1$ in \eqref{eq1}.}

 Under the alternative, the $N_s$  amplitudes $\alpha_q\in \mathbb{R}^{*+}$, frequencies  $f_q\in \mathbb{R}^{*+}$ and phases $\varphi_q \in [0,2\pi[$ of the deterministic part are unknown.  In RV exoplanet detection, this deterministic part represents the planetary signature(s). 
 
%%%%%%%%%%%%%%%
{
Note that model \eqref{hyp} can actually be useful for the detection of periodic signals more general than only pure sinusoids.
In such cases, the Fourier spectrum may contain many harmonics. Because the fundamental frequency has zero probability to coincide with a Fourier frequency, the corresponding number of nonzero Fourier coefficients (\textit{i.e.}, of deviations under the alternative) always equals $N$.
However, most of the energy of periodic signals is captured by a small fraction of Fourier coefficients, so that model \eqref{hyp}  would often be accurate for such signals
 with some $N_s\ll N$. \\
 In the case of RV signals for instance, a study of the influence of the Keplerian parameters (planets' { orbital} parameters) shows that this  is indeed the case \cite{Sulis_2016c}.
% mainly four of them influence the sparsity of the Fourier spectrum \cite{Sulis_2016c}.
% \begin{itemize}
%\item	The eccentricity parameter causes the major deviation to pure sinusoidal signals, as non circular orbits tend to create peaky signals as $e$  approaches one. This in turns translates into largest harmonics in the Fourier space. High-eccentricity orbits are rare, however \cite{Baluev_2015}. In the  Solar System for instance,  the highest planetary eccentricity is $e = 0.058$ for Mercury. 
%
%\item The value of the planet period
%	is related to spectral leakage. If we consider other parameters leading to a sinusoidal signal, all Fourier components
%	will be nonzero unless  the inverse of the period falls exactly on the Fourier frequency grid, which happens with vanishing probability. However, 
%	the energy of the pure sine leaks mainly on a few contiguous frequencies, even when the planet frequency falls half way between two Fourier frequencies.
%
%\item	For high eccentricities, the argument of the periastron $\omega$ modifies the RV curve with, however small impact on the Fourier spectrum.
%
%\item 	The number of planets, $N_p$,  obviously impacts directly the number of components in the Fourier spectrum through the sum in Eq. \eqref{hyp}. This number can however safely be considered in the range $1$ to $10$  for the vast majority of exoplanetary systems.
%\end{itemize}
In all  cases except perhaps very rare and 
exotic systems, the RV spectral signatures exhibit only a small fraction { of} significant harmonic{s} because the planets 
 tend to have low eccentricities and are in small number ($N_p$). In short, the spectrum is sparse (though not strictly sparse) and RV signals can be modeled by a sum of a
small number of  pure sinusoids.
We call this number $N_s$, and we say that $N_s \ll N$. When there is one planet with frequency close to the Fourier grid, $N_s$
will be essentially $1$. In our simulations { for} Sec. VII, we found that  multiplanetary systems with $5$ eccentric planets behave
like model (3) with $N_s$ not exceeding, say,  $20$ at most. 
}
%%%%%%%%%%%%%%%%%
%{\color{blue}
%In fact, the planetary signatures are Keplerian (see \eqref{eq_RV} in Appendix \ref{App0}) but, for low-amplitude exoplanets, we expect to find periodograms with one, or a few harmonic signals. Harmonics will appear in case of off-grid signal frequencies or high eccentricities. But we note that high-eccentricity orbits are rare \cite{Baluev_2015} and if several planets are present, $N_s$ encapsulates all the corresponding harmonic signals.
%The sparsity of the planetary signatures in the Fourier domain will be illustrated in Sec.\ref{SecF}.
%} 

%-----------------------------------------------------------------------------------------------------------------------------------------------
%
\section{Periodograms' statistics: asymptotics}
 \label{sec3}
  
\subsection{Classical (Schuster's) periodogram}
 
The  frequencies considered { in %Eq.
 \eqref{eq1} will be}  Fourier frequencies $\{\nu_k:=\frac{k}{N}\}_{k=0,\hdots,N-1}$ and  $N$ is considered even. For simplicity but without loss of generality we will often consider the subset of $(\frac{N}{2}-1)$  Fourier frequencies corresponding to $k\in \Omega := \{ 1,\hdots,  {\frac{N}{2}-1}\}$.
Asymptotically, the  periodogram $P$ in %Eq.
{ \eqref{eq1}} is an unbiased but inconsistent estimate of the PSD \cite{Brillinger_1981}. %, Chap. 5, theorems 5.2.1 et 5.2.4 
   Under the above assumptions on $E$, the periodogram ordinates at different frequencies $\nu_k$ and $\nu_{k'}$ are asymptotically independent  \cite{Li_2014}. \\
Under  ${\cal{H}}_0$, the asymptotic distribution of $P$ is (Th. 5.2.6,\cite{Brillinger_1981}):
\begin{equation}  P(\nu_k | {\cal{H}}_0) \sim
  \left\{         
      \begin{aligned}
	& \frac{S_E(\nu_k)}{2}  \chi^2_2 , ~~~~ \forall k ~ \in ~ \Omega,   \\
	&  S_E(\nu_k) \chi^2_1, ~ \text{ for } ~ k = 0, \frac{N}{2}.
      \end{aligned}
    \right.
    \label{dist}
\end{equation}
Under ${\cal{H}}_1$, the distribution of $P$ is known when cisoids are present  in (\ref{hyp}) (\cite{Li_2014}, Corollary 6.2(b)).
The real case of  model  (\ref{hyp}) can be treated similarly (see  Appendix \ref{app1}). This leads to
the asymptotic distribution:
\begin{equation}  P(\nu_k | {\cal{H}}_1) \sim
  \left\{         
      \begin{aligned}
	 & \frac{S_E(\nu_k)}{2}  \chi_{2,  \lambda_k}^2 ,~ \forall k ~ \in ~ \Omega,  \\
	 &{S_E(\nu_k)}  \chi^2_{1,  \lambda_k}, ~ \text{ for } ~ k = 0, \frac{N}{2}.
      \end{aligned}
    \right.
    \label{dist2}
\end{equation}
The  $\{\lambda_k := \lambda(\nu_k)\}$ are  non centrality parameters given for $k ~ \in ~ \Omega$ by
\begin{equation} \hspace{-3mm}
  \label{lambda}  
	 \lambda_k \!= \!  \frac{N}{2 S_E(\nu_k)}   \sum_{q=1}^{N_s}  \Big[  \alpha_q^2 \kappa_q ^2 + 2 \alpha_q \kappa_q \sum_{\ell = q+1}^{N_s}  \alpha_\ell \kappa_\ell  \cos(\theta_q-\theta_\ell )\Big],
\end{equation}
and for $k = 0, \frac{N}{2}$ this expression is halved.
The terms $\kappa_q$ and $\theta_q$, given by \eqref{eq_keppa} and \eqref{eq_theta} in Appendix A, arise from  spectral leakage
{  through the spectral window $K_N$  \eqref{eq_Kn}.
Owing to the fast decay of $K_N(\nu)$,
%as a function of $\nu$, 
the proportion of parameters  $\lambda_k$ that significantly differ from $0$ is small if $N_s\ll N$.}
%{
%\\ In case of a Keplerian planetary signatures under ${\cal{H}}_1$, the expression of this non centrality parameter should become extremely more complicated.}

%---------------------------------------------------
\subsection{Averaged periodogram}
We assume that a training data set $\mathcal{T}$ of independent realisations of the { colored }noise is available. This set  is obtained by $L$ independent simulations corresponding to $L$ time series $X_\ell$ sampled on the same grid as the observations:  $\mathcal{T}={\big\{}\{X_\ell(t_j)\}_{j=1,\hdots,N}{\big\}}_{\ell=1,\hdots,L}$.
A straightforward estimate of the noise PSD {  is} the averaged periodogram \cite{Bartlett_1950}:

$$ 
	\overline{P}_L(\nu_k |  {\cal{H}}_0)  : =  \frac{1}{L} \sum_{\ell=1}^{L}  \frac{1}{N} \Big| \sum_{j=1}^{N} X_\ell(t_j)\mathrm{e}^{-{\rm{i}}2\pi\nu_k j} \Big| ^2.
$$

Using \eqref{dist}, the asymptotic distribution of $\overline{P}_L$ {  is:}
\begin{equation}
 \overline{P}_L(\nu_k | {\cal{H}}_0)  \sim
  \left\{         
      \begin{aligned}
	  &\frac{S_E(\nu_k)}{2L}  \chi_{2L}^2,  ~ &\forall k ~ \in ~ \Omega, \\
	   &\frac{S_E(\nu_k)}{L}  \chi_{L}^2, ~ &\text{ for } ~ k = 0, \frac{N}{2}.
      \end{aligned}
    \right.
    \label{dPmoy}
\end{equation}
$\overline{P}_L$ is a consistent and unbiased estimate of $S_E(\nu)$ when both $L\rightarrow \infty$ and $N \rightarrow \infty$ (\cite{Proakis_1996}, Chap.14).

The effect of  stochastic estimation noise caused by the finiteness of $\mathcal{T}$ is encapsulated in $L$. This clearly impacts the distribution
of $\overline{P}_L$ in \eqref{dPmoy} and in turn the efficiency of the subsequent standardization by $\overline{P}_L$. 
%---------------------------------------------------
\subsection{Periodograms standardized with $\overline{P}_L$}
\label{standard}
We now turn to  statistical properties of standardized periodograms of the form \eqref{ptilde}.
When the averaged periodogram $\overline{P}_L$ is used, this yields:
\begin{equation} 
\label{eq_pr}
	 \widetilde{P}(\nu_k \;|\; \overline{P}_L) := \frac{P(\nu_k)}{\overline{P}_L(\nu_k)}.
\end{equation}
\indent As the numerator and denominator are independent variables with known asymptotic distributions, assessing the distribution of their ratio is straightforward. 
The ratio of two  independent random variables (r.v.)  $V_1\!\sim\! \chi_{d_1}^2$ and $V_2\!\sim\! \chi_{d_2}^2$ follows a Fisher-Snedecor law  noted $F(d_1,d_2)$ with ($d_1,d_2$) degrees of freedom: $\frac{V_1/d_1}{V_2/d_2} \sim F(d_1,d_2)$ \cite{Abramowitz_1972}. \\
Consequently, from (\ref{dist}) and (\ref{dPmoy}), the asymptotic distribution of this standardized periodogram under $\mathcal{H}_0$ is:
\begin{equation} \hspace{-3mm}
\widetilde{P}(\nu_k   |  \overline{P}_L, {\cal{H}}_0)  \! \sim \!
  \left\{         
      \begin{aligned}
	   &  \! \frac{S_E(\nu_k)\chi_2^2/2}{S_E(\nu_k)\chi_{2L}^2/2L}  \!  \sim F(2,2L),  \forall k  \in  \Omega ,\\
	    & \!  \frac{S_E(\nu_k)\chi_1^2}{S_E(\nu_k)\chi_{L}^2/L}  \! \sim F(1,L),   \text{ for }  k  \!= \! 0, \frac{N}{2}.
      \end{aligned}
    \right.
    \label{dist3}
\end{equation}
Similarly, from (\ref{dist2}) and  (\ref{dPmoy}), { we have  under $\mathcal{H}_1$}:
\begin{equation}
 \! \widetilde{P}(\nu_k |  \overline{P}_L, {\cal{H}}_1) \! \sim \!
  \left\{         
      \begin{aligned}
	   &\frac{\chi_{2, \lambda_k}^2/2}{\chi_{2L}^2/2L} \sim F_{\lambda_k}(2,2L),   ~ \forall k  \in  \Omega, \\
	    &\frac{\chi_{1, \lambda_k}^2}{\chi_{L}^2/L} \sim F_{\lambda_k}(1,L),  \text{ for }  k = 0, \frac{N}{2},
      \end{aligned}
    \right.
	\label{dist_H1}
\end{equation}
where $F_{\lambda_k}$ denotes a non-central $F$ distribution with non centrality parameter  $\lambda_k$ given by (\ref{lambda}).

Under the null hypothesis,  \eqref{dist3} shows that the distribution of the standardized periodogram is independent of the nuisance signal, \textit{i.e.}, of the partially unknown noise PSD. This  property is  important as it will be inherited by some of the test statistics investigated in Sec. \ref{sec6}, leading to  CFAR tests.
%
%-----------------------------------------------------------------------------------------------------------------------------------------------
\section{Considered tests}
 \label{sec4}
\subsection{Preliminary notations}
{ These tests are better presented using $P$-values\footnote{We will denote the $P$-values 
by $V$ because $P$ denotes the periodogram in this work.} and order statistics. When necessary
the notation will distinguish between a r.v.  $Z$ and its realization $z$. We recall that the observed 
$P$-value $v_Z$ 
%corresponding to the realization $z$ of a random variable $Z$ 
%is
 { is defined as:}
$$
	v_Z(z)  := \textrm{Pr\;} (Z>z)
$$
and $v_Z$ is one realization of the r.v. $V_Z$, which is uniformly distributed.
Similarly,  for a vector of r.v.  ${\bf{Z}}=[Z_1,Z_2,\hdots,Z_N]^\top$ of which ${\bf{z}}=[z_1,z_2,\hdots,z_N]^\top$
 is one realization we  will denote by
 $$ \displaystyle{\min_k\;} {{z_k}}:=z_{(1)}<z_{(2)}< \hdots<{ z_{(N)}} := \max_k \; {{z_k}}$$
 the ordered values of ${\bf{z}}$ and by $Z_{(1)},\hdots,Z_{(N)}$ the  order statistics of  ${\bf{Z}}$.
% $$ \displaystyle{\min_k\;} {{Z_k}}:=Z_{(1)}<Z_{(2)}< \hdots<{ Z_{(N)}} := \max_k \; {{Z_k}}.$$
The observed $P$-values corresponding to ${\bf{z}}$ will be denoted by $v_{{\bf{Z}},k}$ (with $v_{{\bf{Z}},k} := \textrm{Pr\;} (Z_k>z_k)$)
and the observed ordered  $P$-values by $v_{{\bf{Z}},(k)}$. The corresponding r.v. will be denoted  by
 $V_{{\bf{Z}},k}$ and $V_{{\bf{Z}},(k)}$. The ordered $P$-values $V_{{\bf{Z}},(k)}$ are not uniform, because they are order statistics 
 from a uniform distribution, and obviously dependent{\cite{David_2003}}.
%In the sequel, ${\bf{Z}}$ will  be taken as one of the periodograms discussed previously. 

 }
%%%%%%%%%%%
% These tests are better presented using  order statistics, which
%for a vector ${\bf{Z}}=[Z_1,Z_2,\hdots,Z_N]^\top$
% will be denoted by:
% $$ \displaystyle{\min_k\;} {{Z_k}}:=Z_{(1)}<Z_{(2)}< \hdots<{ Z_{(N)}} := \max_k \; {{Z_k}}.$$
%In the sequel, ${\bf{Z}}$ will indeed be taken as one of the periodograms discussed previously. 
%
%We recall that the observed $P$-value\footnote{We will denote the $P$-values by $V$ because $P$ denotes the periodogram in this work.} $v_Z$ corresponding to the realization $z$ of a random variable $Z$ is defined as:
%$$
%	v_Z(z)  := \textrm{Pr\;} (Z>z).
%$$
%The $P$-values $V_Z$ are uniformly distributed. The ordered $P$-values, noted $V_{Z,(k)}$ hereafter,
% are not uniform, because they are order statistics from a uniform distribution, and obviously dependent{\cite{David_2003}}.
%%%%%%%%%%%%%%%

We now present  some  tests discussed in the Introduction and selected for reference as they  cover different { classical models}. 
Under { assumptions specified below} on the distributions of  the variates $\{Z_i\}$,  the properties of these tests are {
 known and we shall summarize them}. 
 %Sec.\ref{sec5}  will describe the counterpart of these tests  when they are applied to $\widetilde{P}$ instead of $P$. 
 %-----------------------------------------------------------------------------------------------------------------------------------------------------------------------------
 \subsection{ { Test} statistics}
 \label{Sec_tests_cl}
 All tests below are of the form
  ${\rm T}({\bf{z}})  \mathop{\gtrless}_{\mathcal{H}_0}^{\mathcal{H}_1} \gamma$,  
  with  { ${\bf{z}}$ the data (which below will be a set of  ordinates of one of the periodograms discussed above),  T$(\cdot)$ the test statistic and $\gamma\in \mathbb{R^+}$  a threshold that determines the false alarm rate. \\
  %Of interested is the distribution of  several   test statistics   $\rm{T}({\bf{Z}}) $. 
}

  \label{sec31}
   \subsubsection{ Test of the maximum}
   \label{Maxtest}

   \begin{equation}
	{\rm {T}_{M}}({\bf{Z}}) :=Z_{(N)}.
 	 \label{lemax}
\end{equation}
For independent variates $Z_i$ of Cumulative Distribution Function (CDF) $\Phi_{Z_i}$, the false alarm of the test is $ {\rm P_{FA}}(\gamma)=1-\prod_{i=1}^N \Phi_{Z_i}(\gamma)$. 
 For a model involving  under $\mathcal{H}_1$ a single sinusoid { with unknown  frequency (but on the Fourier grid) }and under  $\mathcal{H}_0$ a white Gaussian noise { (WGN)} of known variance $\sigma^2$, the ordinates of the periodogram ${\bf{P}}$  evaluated at successive Fourier frequencies are under  $\mathcal{H}_0$   independent and identically distributed (i.i.d.) with distributions given by \eqref{dist}, where $S_E=\sigma^2$.
  For this model  ${\rm T_M}({\bf{P}})$  corresponds to the GLR test  \cite{Kay_1998}. \\
  %This test is powerful when ${\bf{Z}}$ is very sparse, see Sec. \ref{hc} . 
 \subsubsection{Fisher's test}
\begin{equation} 
	\label{eq3b}
	{\rm {T}_{F}}({\bf{Z}}):= \frac{Z_{({N})}}{\displaystyle{\sum_{k ~ } Z_{(k)}}}.
\end{equation}

When applied to periodogram of { WGN} of unknown variance, Fisher's test is CFAR  (the distribution of the test statistics is established in \cite{Fisher_1929}). {This test is actually the GLR test under the model of a single sinusoid on the Fourier grid and WGN of unknown variance (see  \cite{Quinn_1986}, who also covers the case of more than one sinusoid).
%The denominator of (\ref{eq3b}) is an asymptotically unbiased estimate of the PSD (up to a known constant that can be absorbed in the threshold) \cite{Priestley_1981}. So in this case the Fisher's test \eqref{eq3b} can be seen as a standardization of the GLR \eqref{lemax} by the estimated variance.
}
 Examples of works using this test in Astronomy are\cite{Koen_2015a,Koen_2015b,Schwarzenberg_1998,Aittokallio_2001,Guitierrez_2009}. \\
 %------------
  \subsubsection{A test inspired by the tests of Chiu and Shimshoni}
\begin{equation} 
\label{test_ch}  
      \begin{aligned}
{\rm T_C}({\bf{Z}}) := {Z_{({N}-N_C+1)}},
	\end{aligned}
\end{equation} 
%with $N_C\geq 1$ a parameter related to the number of sinusoids.
 {with $N_C\geq 1$ a parameter related to the number of sinusoids. }
  
This test statistic is justified by the observation made in \cite{Shimshoni_1971,Chiu_1989} that for multiple sinusoids, order statistics
different from the maximum may be more discriminative {than $Z_{(N)}$} against the null. 
% { Multiple peaks could be caused by multiple sinusoids of by the leakage from nearby frequencies or harmonics of fundamentals frequencies.}
These tests are designed for periodograms of white noise of unknown variance and their asymptotic false alarm rates are given in  \cite{Shimshoni_1971,Chiu_1989}.
As Fisher's test, they involve denominators whose purpose is normalization by consistent estimates of the noise variance.
${\rm T_C}$ is a simplification of these tests: the normalization is ignored, because it will not be necessary to yield a CFAR detector
once applied to periodograms standardized by $\overline{P}_L$. The false alarm rate of test  ${\rm T_C}$ for white noise can be deduced from the expression obtained in Sec. \ref{TC}.

As for the tests \cite{Shimshoni_1971,Chiu_1989}, we expect ${\rm T_C}$ to have decreasing power in case of strong mismatch between the value of parameter $N_C$ and the  {  number of detectable deviations under   ${\mathcal{H}}_1$ (roughly speaking, $N_s$)}.
Not fixing $N_C$ in advance but estimating this parameter from the data may lead to more powerful tests, but at the cost of a more difficult control of the FA rate (as $N_C$ becomes random). 
This suggests to consider other approaches that are adaptive in the number of sinusoids, which is the case of the last two tests \eqref{HC_test} and \eqref{BJ}.\\

  \subsubsection{Higher Criticism}
  \label{hc}
{This   test statistic  is defined by  \cite{Donoho_2004}}:
{
 \begin{equation}  
		 { \rm HC^{\star}}({\bf{  Z}}):= \!\!\displaystyle{\max_{1 \le { k \le \alpha_0 N}}  \frac{\sqrt{N}(k /N - v_{{\bf{Z}},( k )})}{\sqrt{v_{{\bf{Z}},(k)}(1-v_{{\bf Z},({ k })}})}},
		\label{HC_test}
\end{equation}	
}
where $v_{{\bf Z},(k)}$ are ordered $P$-values and parameter $\alpha_0\in [\frac{1}{N},1]$.

HC is  designed under the assumption that under the null hypothesis the  ordinates $\{Z_k\}$ are  i.i.d. with known marginal distribution.
When ${\bf Z}=\frac{2{ \bf{P}}}{\sigma_E^2}$, with ${\bf P}$ the periodogram of a white noise of known  variance $\sigma_E^2$ under the null, this distribution is given by \eqref{dist} with  $S_E=\sigma_E^2$. 
 {
 The ordered $P$-values involved  in \eqref{HC_test} are thus ordered values of
\begin{equation}
%{
v_{  {2{ \bf{P}}}/{\sigma_E^2}  ,\;k}
% V_{P}(\nu_k)}:=
  \left\{         
      \begin{aligned}
	   &1-\Phi_{\chi^2_{2}}(P(\nu_k)),   ~ \forall k ~ \in ~ \Omega, \\
	    &1-\Phi_{\chi^2_{1}}(P(\nu_k)),   ~ \text{ for } ~ k = 0, \frac{N}{2}.
      \end{aligned}
    \right.
    \label{pval1}    
\end{equation}
}
Under the alternative, a fraction of the ordinates contain a deterministic part, and hence follow \eqref{dist2}, with  $S_E=\sigma_E^2$ ({see Sec. 1.7  of \cite{Donoho_2004}}). The frequencies at which these deviations occur are unknown. Their magnitudes ({ related to} $\lambda_k$) are also unknown and weak, in the sense that they are comparable to the expected magnitude of the periodogram maximum  under the null hypothesis. 

Optimality\footnote{A test is said optimal in \cite{Donoho_2004} if the sum of the probability of error and the probability of missed detection tends to zero.} results of HC are  asymptotic and established for a specific sparsity \textit{vs} amplitude model. 
When the deviation occurs at a single frequency or for extremely sparse signatures in the Fourier domain, 
\cite{Donoho_2004} showed that the test based on the maximum  \eqref{lemax} is asymptotically optimal in the sense that whenever the Neyman-Pearson test has full power, test \eqref{lemax} has full power as well. In such cases, the largest periodogram ordinate (or, equivalently,  the smallest $P$-value) is the most powerful test statistic to discriminate against the null.\\
As discussed about the tests \cite{Shimshoni_1971,Chiu_1989} and ${\rm T_C}$ above, better regions than the maximum ordinate/smallest $P$-value
may be useful for sparse but not extremely sparse signatures. 
\cite{Donoho_2004} showed that there exists sparse signals of weak amplitude that can be  optimally detected (in the sense of full power) by HC but not by ${\rm T_M}$. This makes HC  particularly interesting in the present study.\\
% With  \eqref{HC_test}, HC looks  at ordered $P$-values with increasing order (or equivalently, at periodogram peaks of decreasing magnitude).  
%HC computes  weighted deviations of the empirical CDF of the $P$-values w.r.t. to their theoretical (uniform) CDF under the null. 
%The parameter $\alpha_0$ expresses the fact that computing all ratios in ${\rm HC^\star}$ is not useful since the alternative is sparse. 
%In practice, this parameter can be set to $1/2$ without impact on the test behaviour.
%-----------------------------------------------------------------------------------------------------------------------------------------------------------------------------
 \subsubsection{Berk-Jones }
 \label{bj}
A  test statistic related to $\rm HC$ is that of Berk-Jones \cite{Moscovich_2016,Aldor_2013,Mary_2014,Kaplan_2014,Gontscharuk_2014} defined by:
{
 \begin{equation}  
	{ \rm BJ({\bf{Z}})}:= \!\!\displaystyle{\max_{1 \le k \le \alpha_0 N}} I_{1-v_{{\bf Z},(k)}} (N-k+1,k), 
	\label{BJ}
\end{equation}	
}
where $I$ denotes the regularized incomplete beta function  \cite{Abramowitz_1972}.\\
The deviations in form of $Z-$scores in \eqref{HC_test} for HC are established using  the asymptotic convergence of a binomial distribution to a Gaussian distribution. In the tails, however, this convergence is very slow (see{\cite{Mary_2014, Li_2015,Moscovich_2016}} for illustrations). For this reason the test statistics \eqref{BJ}, which compares favorably to HC and other goodness-of-fit (GOF) tests in some cases, was recently   (and almost simultaneously) proposed by{ \cite{Moscovich_2016,Aldor_2013,Kaplan_2014,Gontscharuk_2014}}. As noted in{ \cite{Moscovich_2016, Gontscharuk_2016}},  this test was initially proposed by Berk and Jones (and called $M_n^+$) in \cite{Berk_1979}.\\
This test is based on the exact significance reflected by the $P$-values, that is, on the $P$-values of the ordered $P$-values. Since the ordered $P$-values are Beta distributed, with  { $V_{{  {2{ \bf{P}}}/{\sigma_E^2} },(k)}\sim \text{Beta}(k,N-k+1)$, their $P$-values involve the CDF of Beta variables, which is an incomplete Beta function: $\Pr(\text{Beta}(k,N-k+1)\leq x)=I_x(k,N-k+1)$ } and leads to test \eqref{BJ} with the $P$-values computed as in \eqref{pval1}.\\
 $\rm BJ$ presents the same adaptive optimality as $\rm HC$ for sparse mixture detection, and the asymptotic distribution of the $\rm BJ$ test statistic  can be found in Th. 4.1 of  \cite{Moscovich_2016}.  As for $\rm HC$, convergence to the asymptotic distribution may be slow. { Efficient algorithms for computing significance levels  (and hence the function $\gamma \mapsto {\rm P_{FA}}(\gamma)$)  of $\rm HC$ and $\rm BJ$ for finite (but possibly large) values of $N$ can be found in\cite{Moscovich_2016,Moscovich_2015b}. $\rm BJ$ and $\rm HC$  are studied from the viewpoint of  local levels in {\cite{Gontscharuk_2016}}.}
%--------------------------------------------------------%
\section{ Tests { applied to $\bf{Z}=$}  $\widetilde{\bf P}\; |\; \overline{\bf P}_L$}
 \label{sec5}

To simplify the presentation of the results we restrict in the following to the frequency set  $\Omega$, \textit{i.e.},
to standardized vectors
$$
{\bf{\widetilde{P}\;|\;\overline{P}}}_L:=\left[\frac{P(\nu_1)}{\overline{P}_L(\nu_1)},\hdots, \frac{P(\nu_{\frac{N}{2}-1})}{\overline{P}_L(\nu_{\frac{N}{2}-1})}\right]^\top.
$$
Extension of the results
to { $\nu_0,\nu_{\frac{N}{2}}$} can be obtained using the distributions \eqref{dist3} and \eqref{dist_H1}.

{
In the following we evaluate false alarm and detection rates by postulating independence of the considered ordinates.
While the asymptotic independence of the periodograms ordinates at positive Fourier frequencies is well established
(\textit{e.g.} \cite{Fan_2003}, Th. 2.14), theoretical results regarding the joint distribution of periodogram ordinates under departures from whiteness (and Gaussian) assumptions are lacking, however. For some results on the largest order statistics, see \cite{Turkman_1984}, who consider MA Gaussian processes, and \cite{Davis_1999}, who consider non-Gaussian sequences. For finite values of $N$, the approaches of \cite{Rife_1974,Quinn_1994} 
might be followed to better characterize the performances of some tests considered below.

In practical situations, the marginal distributions of the considered ordinates are only approximated by their asymptotic distribution 
(this is visible in Eq. \eqref{lerho} for instance).
Besides, as $N$ grows the correlation between the periodogram ordinates at the signal frequencies approaches zero at a slower rate than the correlation between other ordinates ($\!\!$\cite{Li_2014}, Theorem 6.5), indicating that departure to the independence assumption may be more pronounced under ${\mathcal{H}}_1$ than ${\mathcal{H}}_0$. Consequently,  we are using the independence as an operational assumption to quantify  the tests' performances, and the validity of the resulting expressions  below should be checked against numerical simulations. Sec. \ref{sec7} provides several examples.   }
 %%%%%%%%
 \subsection{${\rm T_M}$.}
Under ${\cal{H}}_0$ and in the asymptotic regime, ${\rm T_M}({\bf \widetilde{\bf P}}\;|\;\overline{\bf P}_L)$ 
is the maximum of independent variables given by \eqref{dist3}. For $k \in \Omega$ these $\frac{N}{2}-1$ variables follow an $F(2,2L)$ distribution with general density given in \cite{Abramowitz_1972}. Using the beta function ${\cal{B}}(1,L):=\int_0^1(1-t)^{L-1}\textrm{d}t=\frac{1}{L}$, this density $\varphi_F(\gamma,2,2L)$ can be expressed as
$$
\varphi_F(\gamma,2,2L)=\frac{1}{{\cal{B}}(1,L)}\cdot\frac{1}{L}\cdot \Big(1+\frac{\gamma}{L}\Big)^{-L-1}=\Big(1+\frac{\gamma}{L}\Big)^{-L-1}.
$$
It can be checked that $\int_0^\infty\varphi_F(\gamma,2,2L)\mathrm{d}\gamma=1$. 
The corresponding CDF  $\Phi_F(\gamma,2,2L)$ is obtained  by  integration of $\varphi_F$:
\begin{equation} 
	\Phi_F(\gamma,2,2L) = \displaystyle \int_{0}^{\gamma} \varphi_F(\gamma,2,2L) \mathrm{d}\gamma = 1 - \Bigg(\frac{L}{L+\gamma} \Bigg)^L.
	\label{laphi}
 \end{equation}
The probability of false alarm ($\rm{P_{FA}}$) can be computed thanks to the  asymptotic independence of the ordinates of ${\widetilde{\bf P}\; |\; \overline{\bf P}_L}$:
\begin{equation} 
\begin{aligned} 
	&\!\rm{P_{FA}}({{\rm T_M}} ({\bf \widetilde{\bf P}}\;|\;\overline{\bf P}_L) ,\gamma) := \textrm{Pr}\; ({\rm T_M}({\bf \widetilde{\bf P}}\;|\;\;\overline{\bf P}_L)\! > \!\gamma | {\cal{H}}_0\!) \\
		&=1 -\displaystyle{ \prod_{k\in \Omega}} \textrm {Pr}\, \left(\widetilde{P}({ \nu_k})  \leq \gamma | {\cal{H}}_0,\overline{P}_L\right) = 1 - \Big(\Phi_F(\gamma,2,2L)\Big)^{\frac{N}{2}-1}\\
%		&= 1 - \Big(\Phi_F(\gamma,2,2L)\Big)^{\frac{N}{2}-1} \\
		&= 1-  \Big( 1-\Big(\frac{L}{\gamma+L}\Big)^L\Big)^{\frac{N}{2}-1}.
\end{aligned}
\label{pfa}
 \end{equation}
 
The probability of false alarm  is { (asymptotically)} independent of the  noise PSD, which makes   ${\rm T_M}({\bf \widetilde{\bf P}}\;|\;\overline{\bf P}_L)$ a CFAR detector. 

 Under ${\mathcal{H}_1}$, using (\ref{dist_H1}) and the approximate independence of periodogram ordinates (\cite{Li_2014}, Theorem 6.5), the probability of detection ($\rm{P_{DET}}$) of ${\rm T_M} ({\bf \widetilde{\bf P}}\;|\;\overline{\bf P}_L)$ can be  approximated as: 
  \begin{equation} 
  \begin{aligned} 
	\!\!\!\!\rm{P_{DET}}({{\rm T_M}} ({\bf \widetilde{\bf P}}\;|\;\overline{\bf P}_L) ,\gamma)&:=  \textrm{Pr} \left({\rm T_M}({\bf \widetilde{\bf P}}\;|\;\;\overline{\bf P}_L)\! > \!\gamma | {\cal{H}}_1\! \right) \\
	&\approx 1 - \displaystyle{ \prod_{k \in \Omega}} \Phi_{F_{\lambda_k}}(\gamma, 2,2L).
	\end{aligned}
		\label{Eq_Pdet}
\end{equation}	
With (\ref{pfa}),  the relationship $\gamma(\rm{P_{FA}})$ for ${{\rm T_M}}$ can be derived as:
 \begin{equation} 
\label{Eq_gam}
		\gamma({{\rm T_M} ({\bf \widetilde{\bf P}}\;|\;\overline{\bf P}_L}),\rm{P_{FA}}) =  L  \Big[ \Big( 1 - ( 1 - \rm{P_{FA}})^{\frac{1}{\eta}} \Big)^{-\frac{1}{L}} -1 \Big],   
\end{equation}	
{where $\eta := \frac{N}{2}-1$}. With (\ref{Eq_Pdet}) and (\ref{Eq_gam}), we  deduce $\rm{P_{DET}(P_{FA})}$ which can be used to compute ROC (Receiver Operating Characteristics) curves:
 \begin{equation} 
   \begin{aligned} 
\label{Eq20}
	&{\rm{P_{DET}}({{\rm T_M}} ({\bf \widetilde{\bf P}}\;|\;\overline{\bf P}_L)},{\rm{P_{FA}}})  \\
	&\approx 1 - \displaystyle{ \prod_{k\in \Omega }}  \Phi_{F_{\lambda_k}}(L  [ ( 1 - ( 1 - \rm{P_{FA}})^{\frac{1}{\eta}} )^{-\frac{1}{L}} -1 ],2,2L)
		\end{aligned}
\end{equation}	
 %-----------------------------------------------------------------------------------------------------------------------------------------------------------------------------
 \subsection{${\rm T_F}$.}
 \label{fisher}
%The Fisher's test ${\rm T_F}$ applied to  ${\bf \widetilde{\bf P}}\;|\;\overline{\bf P}_L$ looks under $\mathcal{H}_0$ for the maximum of identically distributed variables, which from  (\ref{dist3}) and (\ref{eq3b}), 
{ Under $\mathcal{H}_0$,} the Fisher's test ${\rm T_F}$ applied to  ${\bf \widetilde{\bf P}}\;|\;\overline{\bf P}_L$ looks for the maximum of identically distributed variables, which from  (\ref{dist3}) and (\ref{eq3b})
%are distributed as:
%$$
%	\begin{aligned}
%	 \frac{\frac{{P}(\nu_i)}{\overline{P}_L(\nu_i)}}{\displaystyle{\sum_{k}\frac{P(\nu_k)}{\overline{ P}_L(\nu_k)}}}  \stackrel{}{\sim} \frac{F(2,2L)}{  \displaystyle{\sum_{k }} {F(2,2L) }}. \\	
%	\end{aligned}
%$$
%
%Each such variable is 
{ correspond to the} ratio  of a $F$ variable over a sum of $F$ variables. To our knowledge, there is no { analytical} characterization of the resulting distribution
{ for finite values of $L$. }
%When $L \to \infty$, as if $X \sim F(d_1,d_2)$ then $\displaystyle{\lim_{d_2 \to \infty}} d_1 X \sim \chi_{d_1}^2$, the $\rm P_{FA}$ derived by Fisher \cite{Fisher_1929} becomes true.} \\
%
Hence, although this test is CFAR, { computing} the false alarm rate is problematic. Resorting to Monte Carlo simulations to evaluate the function $\gamma \mapsto \rm{{P}_{FA}}$ is not possible, owing to the limited number of available noise realizations.
 Similar remarks can be made about standardized versions of other tests like \cite{Shimshoni_1971, Siegel_1980,Chiu_1989}.
%-----------------------------------------------------------------------------------------------------------------------------------------------------------------------------
\subsection{${\rm T_C}$.}
\label{TC}
We turn to ${\rm T_C}({\bf \widetilde{\bf P}}\;|\;\overline{\bf P}_L,N_C) $, where parameter $N_C$ allows to focus on the regions of order statistics 
where  deviations under the alternative are likely to be  significant. The $\rm{P_{FA}}$ can be obtained by observing
 that, owing to \eqref{dist3}, the number $K$ of standardized ordinates larger than $\gamma$  under $\mathcal{H}_0$ follows a binomial distribution:
 {
% \begin{equation*}
$	 K\sim \rm{Bin}(N/2-1,1-\Phi_{F}(\gamma,2,2L)),$
% \end{equation*}
  whose CDF is: %a $ I$ function:
% $$
 $	\textrm{Pr}\;(K\leq k)=I_{\Phi_{F}(\gamma,2,2L)}(N/2-N_C,N_C), $ see \cite{Abramowitz_1972}.
%$$  
}
 Using \eqref{test_ch}, (\ref{laphi}) and noting
\begin{equation*}  
 	u :=  1-\Phi_{F}(\gamma,2,2L) = \Big( \frac{L}{\gamma+L} \Big)^L,
 \end{equation*}
the $\rm{P_{FA}}$ of this test applied to ${\bf \widetilde{\bf P}}\;|\;\overline{\bf P}_L$ can be expressed as:
\begin{equation}  
\begin{aligned} 
	&\!\!\rm{P_{FA}}({\rm T_C}( {\bf \widetilde{\bf P}}|\overline{\bf P}_L, N_C),\gamma) \!:=\! \textrm{Pr}\; ( {T_C  ({\bf \widetilde{\bf P}}|\overline{\bf P}_L,N_C)}\!\! > \!\gamma | {\cal{H}}_0\!)\! \\
	&= 1- { \sum_{k=0}^{N_C-1} } \textrm{Pr}\; ( K=k\;|\;{\cal{H}}_0) = 1 -  I_{1-u}( \frac{N}{2}-N_C, N_C)  \\
	%&= 1 -  I_{1-u}( \frac{N}{2}-N_C, N_C) \\
	&=   I_{u} (N_C, \frac{N}{2}-N_C),
\end{aligned}
\label{pfa_Tor}
 \end{equation}
 where the last equation uses   $I_x(a,b)=1-I_{1-x}(b,a)$  (\textit{cf} Prop. 6.6.3 in \cite{Abramowitz_1972}).
%-----------------------------------------------------------------------------------------------------------------------------------------------------------------------------
%
As ${\rm T_M}({\bf \widetilde{\bf P}}\;|\;\overline{\bf P}_L)$, this test is CFAR.

As a side remark, note that the false alarm rate of the test ${\rm T_C}$ applied to the periodogram of a white noise of known variance is obtained as above,
by replacing $\Phi_F$ with the CDF of $\chi^2$ variables according to \eqref{dist}.

 The function $\gamma \mapsto \rm{P_{DET}} (\gamma)$ of ${\rm T_C}( {\bf \widetilde{\bf P}}\;|\;\overline{\bf P}_L,N_C)$  can be deduced similarly to (\ref{pfa_Tor}), with the difference that $K$ is no longer binomial owing to the $\{\lambda_k\}$.  Appendix \ref{app6} shows that
 \begin{equation}  
 \begin{aligned} 
 		&\rm{P_{DET}}( T_C( {\bf \widetilde{\bf P}}\;|\;\overline{\bf P}_L,N_C),\!\gamma):= \textrm{Pr}\; (\!  T_C( {\bf \widetilde{\bf P}}\;|\;\overline{\bf P}_L,N_C)\!  > \gamma | {\cal{H}}_1) \\
		&\approx 1 -\!\! \displaystyle{ \sum_{i=0}^{N_C - 1} \sum_{\Omega^{(i)} \in \Omega^i} \prod_{k=1}^{ i  } } \Big(1 - \varphi_{F_{\lambda_{\Omega^{(i)}_k}}\!^{\!\!\!\!\!\!\!\!\!(\gamma, 2,2L)}} \Big) \displaystyle{ \!\!\!\! \prod_{k'=1}^{ \frac{N}{2}-1-i }}\!\!\varphi_{F_{\lambda_{ \overline{\Omega}^{(i)}\!_{\!\!\!\!\!\!k'} }}{(\gamma, 2,2L)} },
		\label{lagrosse}
\end{aligned} 
\end{equation}	
which can be used with \eqref{pfa_Tor} to compute ROC curves. The non centrality parameters $\{\lambda_{\Omega_k^{(i)}}:= \lambda(\nu_{\Omega_k^{(i)}})\}$ and $\{\lambda_{\overline{\Omega}_{k'}^{(i)}}:= \lambda(\nu_{\overline{\Omega}_{k'}^{(i)}})\}$ are given by (\ref{lambda}) with the notation defined in \eqref{pfac}.
{ Note that the tests $\rm T_M$ and $\rm T_C$  are both working on order statistics of  $P$-values, which are Beta random variables, so in both cases the $\rm{P_{FA}}$ is given by the tail of a Beta distribution. }

%-----------------------------------------------------------------------------------------------------------------------------------------------
\subsection{${\rm HC^\star}$ and ${\rm BJ}$.}
 \label{sec_HC}
 When applied to ${\bf \widetilde{\bf P}}\;|\;\overline{\bf P}_L$, the $P$-values involved in the HC test \eqref{HC_test} and in the BJ test 
 \eqref{BJ} should be computed according to the distribution under the null given by \eqref{dist3}. Hence,  \eqref{pval1} is replaced by:
\vspace{-0.4cm}
{
$$
\vspace{-0.2cm}
{ v_{{\bf \widetilde{\bf P}}|\overline{\bf P}_L,\;k}}  :=
  \left\{         
      \begin{aligned}
	   &1-\Phi_F\left(\frac{P(\nu_k)}{\overline{P}_L(\nu_k)}, 2, 2L\right) ,   ~ \forall k ~ \in ~ \Omega, \\
	    &1-\Phi_F\left(\frac{P(\nu_k)}{\overline{P}_L(\nu_k)}, 1, L\right),   ~ \text{ for } ~ k = 0, \frac{N}{2}.
      \end{aligned}
    \right.
$$
}
The properties of the two tests are otherwise left unchanged, with the CFAR property added:  thanks to the standardization of ${{\bf P}}$ by $\overline{\bf P}_L$, the $P$-values are
independent of the noise PSD. 
%{ Indiquer ici comment calculer ces p-values, cf. com referee.}

%%%%%%%%%

\section{ Tests { applied to ${\bf{Z}}=$ }  $\widetilde{\bf P} \;|\;\widehat{\bf S}_E$}
\label{sec6}
Estimates of  $\widehat{S}_E$ different from the  averaged periodogram $\overline{P}_L$ can be used for standardization in \eqref{eq_pr}. 
Sec. \ref{rel} has reviewed some  methods to obtain such estimates.
For the purpose of comparing  $\widetilde{\bf P} |\widehat{\bf S}_E$ with  ${\bf \widetilde{\bf P}}|\overline{\bf P}_L$,   
we opt for parametric estimates  allowing { for} an automatic parameter setting, as this approach is  commonly used in practice.
For instance, when using an estimate based on an  AR process of order $o$, this order can be estimated using many criteria, {\textit{e.g.}} \cite{Akaike_1969, Akaike_1974,Parzen_1975, Hannan_1979, Rissanen_1984}. 

In our studies, we found that the selected order is often different from the true order (as expected, see \textit{e.g.} { p. 211 of \cite{Akaike_1969} and \cite{Boardman_2002}} for similar conclusions) but, as far as detection results are concerned,  these criteria  have very similar  behaviour for sufficiently large $N$.   

In any such method, let us denote respectively  by $\widehat{o}_{AR}$, $\{\widehat{c}_j\}_{j=,1,\cdots,\widehat{o}_{AR}}$ and 
$\widehat{\sigma}^2(\widehat{o}_{AR})$ the selected order,  AR coefficients and corresponding estimated prediction error variance.
{ The  PSD estimate and resulting standardized periodogram are% is then
 \begin{equation}  \vspace{-0.2cm}
	\widehat{S}_{E,AR}(\nu)\!: =\! \frac{\widehat{\sigma}^2(\widehat{o}_{AR})}{\Big|1\!+\!\displaystyle{\sum_{j=1}^{\widehat{o}_{AR}}} \widehat{c}_j \mathrm{e}^{-2\pi i j \nu }\! \Big|^2},
%\end{equation}	
%and the standardized periodogram  \eqref{ptilde} becomes
%\begin{equation} 
%\label{eq_pr2}
	 \widetilde{P}(\nu_k | \widehat{S}_{E,AR}  )\! :=\! \frac{P(\nu_k)}{\widehat{S}_{E,AR} (\nu_k)}.
	 \label{eq_SAR}
\end{equation}
}

Even if such approaches are relatively straightforward to implement, characterizing the distribution of $ \widetilde{P}(\nu_k |\widehat{S}_{E,AR}  )$
 is  more difficult than in the case of  $ \widetilde{P}(\nu_k | \overline{P}_L  )$, owing essentially to the stochastic nature of $\widehat{o}_{AR}$ in \eqref{eq_SAR}. 
In practice, the `whitening' effect of  $\widehat{S}_{E,AR}$ is efficient because the selection procedures have good fitting properties (they are approximately consistent,  \textit{i.e.},  $\hat{S}_{E,AR} \overset{\approx}{\longrightarrow}S_E$ as $N \longrightarrow+\infty$).
 This leads to consider as reasonable  the assumptions that, in effect,  $ \frac{{ { P (\nu_k)}}}{\widehat{ S}_{E,AR} (\nu_k)}\approx  \frac{P(\nu_k)}{S_E(\nu_k)}$
 and, with \eqref{dist}, that $\frac{P(\nu_k)}{S_E(\nu_k)}$ is approximately a  $\chi_2^2/2$ r.v. for $  k \in \Omega$ and a ${\chi_1^2}$ r.v. for $  k= 0,\frac{N}{2}$.
An approximate false alarm rate can then be evaluated from these assumptions.
For example, for the ${\rm T_C}$ test applied to $\frac{\bf{P}}{\widehat{ \bf{S}}_{E,AR} }$, following the lines of \eqref{pfa_Tor} leads to:

\begin{equation}\vspace{-0.1cm}
	\begin{aligned} 
	\label{PFA_param}
 	 &{\rm{P_{FA}}}({\rm T_C}( {\bf \widetilde{\bf P}}|{\widehat{ \bf S}_{E,AR}},N_C ),\gamma)  \approx  {\rm{P_{FA}}}({\rm T_C}( {\bf \widetilde{\bf P}}|{\bf S}_E, N_C),\gamma) \\
%	 &~~~~~~~~~~~~~~~~~~~~~~~~~ \approx  {\rm{P_{FA}}}({\rm T_C}( {\bf \widetilde{\bf P}}|S_E, N_C),\gamma) \\
 	 &~~~~~~~~~~~~~~~~~~~~~~~~~~~~~~~~ = I_{u} (N_C, \frac{N}{2}-N_C) .
 	 \end{aligned}
\end{equation}

with $ u := \Phi_{\chi_2^2}(2\gamma)= \mathrm{e}^{-\gamma}$.
 We will evaluate the reliability of this approximation by numerical simulations in Sec.\ref{SecD}.
%{ We will also show that, even if the parametric model is the true one (e.g; AR colored noise), the resulting $\rm{P_{FA}}$control is less reliable than the one derived by non parametric approach even if the tests' performances can be more powerful (e.g; for low $L$ value). }
%-----------------------------------------------------------------------------------------------------------------------------------------------
\section{Numerical study}
 \label{sec7} 

 \subsection{Simulation setting}
\label{SecA}
%We use here data generated according to model \eqref{hyp} in order to validate the theoretical analyses and to compare the performances of the different detection methods discussed in Sec. \ref{sec4} to \ref{sec6}.
%The use of a synthetic noise is necessary  to perform Monte Carlo (MC) simulations, as currently available  time series of simulated stellar noise   are in too small number for that purpose 
{ 
%and real RV observations to detect exoplanet are irregular. 
We consider under $ {\cal{H}}_0$ two PSD models for the noise $E$.
%
%\begin{enumerate}
%\item 
The first model comes from real RV data of the Sun  
%{ observed solar RV data} 
obtained from { the} GOLF spectrophotometer on board SoHO satellite \cite{Garcia_2005}. This instrument has been observing the Sun for 18 years with a sampling rate of 20 s. As several gaps are present in the resulting time series, we selected some ($158$) regularly sampled data blocks of $T \approx 23$ days, of which we averaged and smoothed the periodograms (see Fig. \ref{Fig0},  left panel).
% of shows the averaged periodogram (blue) of the $158$ available blocks (one of these blocks is also represented in grey) sampled each $30$ minutes (to avoid the determinist part of the solar oscillations at high frequencies). The resulting PSD has been filtered using a moving average filter on [3, 10, 20] frequency values respectively for the frequency bands : [0, 2e-5];  [2e-5, 5e-5]; and  [5e-5, 3e-4] Hz.
	(Note that the data we used are filtered at low frequencies so that the resulting PSD estimate does probably not accurately reflect the solar { PSD} at low frequencies).\\
	 The second model (Fig. \ref{Fig0},  right panel) corresponds to a zero-mean second-order stationary Gaussian AR(6) process.
	 % defined by:
	%	\begin{equation*} 
   	%		   \begin{aligned}
	%				 E(t_k) & = 0.7\, E({t_{k-1}}) +  0.05 \,E({t_{k-2}}) + 0.3\, E({t_{k-4}} )  \\
	%				 &- 0.3\, E({t_{k-6}} ) + W({t_{k}}),
     	%		 \end{aligned}
	%	\end{equation*}
%where $t_k:=k\,\Delta t$ and $W$ is a WGN with zero mean and variance $\sigma_W^2\! = \!1$.
%The AR coefficients have been arbitrary chosen and respect the stationary conditions ($\in ]-1,1]$). 
%The noise PSD of the considered AR process (shown in log scale in Fig. \ref{Fig0}, right panel) 
The coefficients were chosen to yield a correlated process exhibiting higher energy at low frequencies  and local variations, as in some stars,
% In the scope of the present study, this strongly colored noise is intended to  roughly represent the large variations encountered  in some stellar PSDs, 
but the choice of this { PSD} is  not intended to reflect  the reality of a particular star. 
%\end{enumerate}
  \begin{figure}[htb!]     \centering
	\includegraphics[width= 9cm,height=4.5cm]{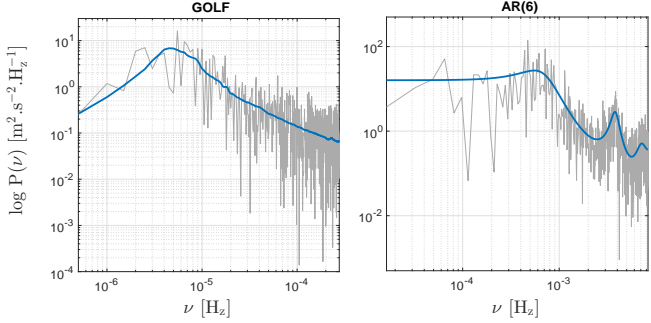}
	\caption{  Left: Estimated PSD of the solar noise with part of GOLF data (blue) and periodogram of one of the $158$ data blocks (grey). Parameters: $N = 1110$, $\Delta t = 30$ min. This PSD  is used to generate the noise for Figs. \ref{Fig2} and \ref{Fig5}. Right: Theoretical PSD of the AR(6) noise (blue) and one noise periodogram (grey). This PSD is used for Fig. \ref{Fig3}. Parameters:$N = 1024$, $\Delta t = 1$ min.}
	\label{Fig0}
		\vspace{-0.2cm}
\end{figure}

}
%

%Under ${\cal{H}}_1$, we consider :the most general case of multiple sinusoids ($N_s \!>\! 1$) whose frequencies $\{f_q\}$ do not belong to the Fourier frequency grid. Several cases will be considered for signal frequencies and amplitudes in the following experiments.  The phase terms $\{\varphi_q\}$ are drawn from a uniform distribution over {$[0, 2\pi[$.} 
{  Under ${\cal{H}}_1$, several cases of exoplanetary RV signatures will be considered.
% for signal frequencies and amplitudes arising from in the following experiments
%% we consider a model evolving as \eqref{hyp}:
%\begin{itemize}
%	\item the most general case of multiple sinusoids ($N_s \!>\! 1$) whose frequencies $\{f_q\}$ do not belong to the Fourier frequency grid (Fig.\ref{Fig2}-\ref{Fig3}). Several cases will be considered for signal frequencies and amplitudes in the following experiments.  The phase terms $\{\varphi_q\}$ are drawn from a uniform distribution over {$[0, 2\pi[$.} 
%	\item the particular case of Keplerian planetary signatures (Fig.\ref{Fig5}-\ref{Fig6}) is described in Appendix \ref{App0}.
%\end{itemize} 
}
The tests' performances will be illustrated by ROC curves representing  $\rm P_{DET}$ as a function of $\rm P_{FA}$.
%For MC simulations, the time series have { $N =1110$ and $1024$ points respectively for the two PSD} and unless otherwise stated $10^4$  realizations have been used. 
{   For { Monte Carlo (MC)} simulations,  $10^4$  realizations have been used. }
%------------------------------------------------------
  \subsection{Analytical expressions for tests based on ${\bf \widetilde{\bf P}}\;|\;\overline{\bf P}_L$ }
  \label{SecB}
\begin{figure}[htb!]   
%\hspace{-0.5cm}
\centerline{\includegraphics[width= 9.5cm,height=9cm]{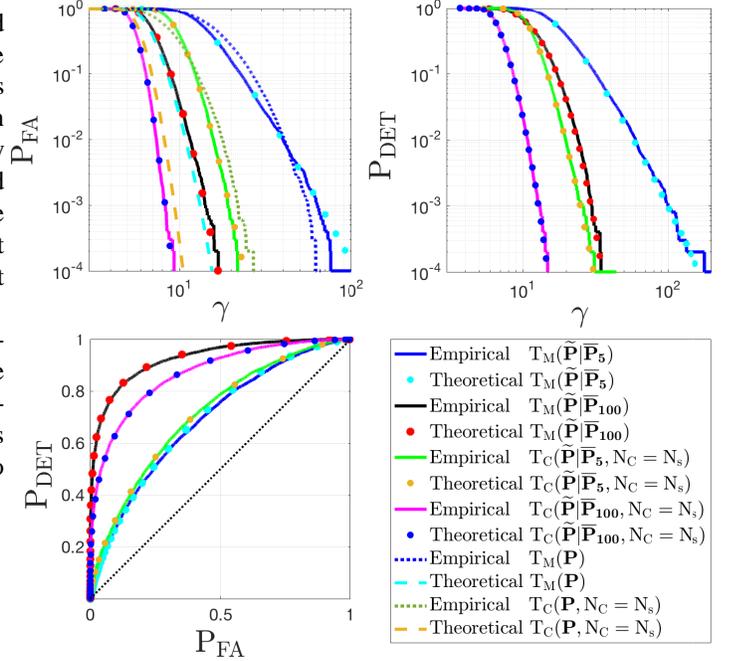}	}
	\caption{ Theoretical \textit{vs} empirical results for ${\rm T_M}( {\bf \widetilde{\bf P}}\,|\,\overline{\bf P}_L)$ and ${\rm T_C}( {\bf \widetilde{\bf P}}\;|\;\overline{\bf P}_L)$ by MC simulations. { Parameters: $\Delta t = 30$ min, $N_s =  3$, $\alpha_q = 0.2$ m.s$^{-1}$ for the three RV signatures of respective periods $11$ h, $2.45$ d and $6.61$ d. The curves in dashed lines show the inconsistency of ignoring noise correlations (Sec. \ref{secC}). }}
		\vspace{-0.4cm}
	\label{Fig2}
\end{figure}

We first consider the tests ${\rm T_M}( {\bf \widetilde{\bf P}}\,|\,\overline{\bf P}_L)$ and ${\rm T_{C}}( {\bf \widetilde{\bf P}}\;|\;\overline{\bf P}_L)$.  
%The  theoretical expressions $\rm{P_{FA}}(\gamma)$, $\rm{P_{DET}}(\gamma)$ and $\rm{P_{DET}}(\rm{P_{FA}})$ corresponding to (\ref{pfa}),  (\ref{Eq_Pdet}), (\ref{Eq20}), (\ref{pfa_Tor}) and (\ref{lagrosse})  are shown with colored dots  in Fig. \ref{Fig2} for $N_s = 3$  sinusoidal signals.  The empirical results from MC simulations are plotted in full lines. Different values of $L$ are considered to illustrate the improvement brought by larger training data sets.
%
{ The first panel of Fig.\ref{Fig2} compares to empirical results obtained from MC simulations to the expressions obtained for the  $\rm P_{FA}(\gamma)$ for both tests (see (\ref{pfa}) and (\ref{pfa_Tor})). The second panel regards the corresponding expressions for $\rm P_{DET}(\gamma)$    (\eqref{Eq_Pdet} and  \eqref{lagrosse}) and the last panel the  expression   \eqref{Eq20} for $\rm P_{DET}(P_{FA})$. All the theoretical expressions are shown by color dots and the empirical results from MC simulations are plotted in full lines. }
Different values of $L$ are considered to illustrate the improvement brought by larger training data sets.
The figure shows a fair agreement between theoretical and empirical results, even for the { not} so large value of $N$ considered here { (N = 1110)}. The test performances logically increase with $L$, as the estimation noise decreases with the increasing size of the training data set. 
%-----------------------------------------------------
\subsection{Effects of standardization by $\overline{\bf P}_L$}
 \label{secC}
We compare the detection performances of the tests ${\rm T_M}$ and ${\rm T_C}$ { applied to} the standardized periodogram $(\widetilde{\bf P}|\overline{\bf P}_L)$ with their unstandardized versions ${\rm T_M}(\bf P)$ and ${\rm T_{C}}(\bf P)$, as described in Sec. \ref{Sec_tests_cl}.
We evaluate first how accurate would be the false alarm rate obtained by neglecting noise correlation for tests ${\rm T_C}$ and ${\rm T_M}$.
For this we assume the detectors consider the noise is white and have knowledge of $\sigma^2_E$, the exact variance of $E$. When ${\bf{Z}}=2 {\bf{P}}/\sigma^2_E$, it is easy to show  that the  false alarm rates assumed by these two tests are  given by 
\begin{equation}
	\!\!\!\!
	\left\{
	\begin{aligned}
		&{\rm P_{FA}}\left({\rm T_M}(2 {\bf{P}}/\sigma^2_E),\gamma\right)=1- \Phi_{\chi^2_2}^{\frac{N}{2}-1}(\gamma),\\
		&{\rm P_{FA}}\left({\rm T_C}(2 {\bf{P}}/\sigma^2_E, N_C),\gamma\right)= I_{\Phi_{\chi^2_2(\gamma)}} (N_C, \frac{N}{2}-N_C).
	\end{aligned}
	\right.
	\label{vs}
\end{equation}
These expressions are compared to the true false alarm rates  in the top left panel of Fig. \ref{Fig2}.
The cyan { dashed} and blue  { dotted} curves show respectively approximation  \eqref{vs} and empirical $\rm P_{FA}$ for ${\rm T_M}$,
 while the yellow  { dashed} and red  { dotted} curves show respectively approximation  \eqref{vs} and empirical $\rm P_{FA}$ for ${\rm T_C}$.
 Clearly, the correspondence between the thresholds values and the target false alarm rates
is destroyed because of noise correlation in absence of standardization. 
%Another remarkable effect, not illustrated here for lack of space, is that false alarms for unstandardized tests accumulate precisely in the regions where the noise PSD is large. Conversely,
%such tests have vanishing detection power in the PSD `valleys'. In contrast, because of standardization by $\overline{\bf P}_L$, 
 %the tests ${\rm T_M}(\widetilde{\bf P}|\overline{\bf P}_L)$ and ${\rm T_{C}}(\widetilde{\bf P}|\overline{\bf P}_L)$  have constant false alarm rate over the whole frequency band{\cite{Sulis_2016a}}. 
%------------------------------------------------------
\subsection{Effects of standardization by $\widehat{\bf S}_{E,AR}$}
 \label{SecD}
As discussed in Sec.\ref{sec6}, a way to deal with the frequency dependence of the noise is to estimate its PSD by parametric models. 
A difficulty with such methods is the injection of estimation noise in the { detection}  process. A standard approach discussed in Sec.\ref{sec6} is to consider that the estimates are sufficiently accurate for their {intrinsic error} to be negligible. We study the performances of this approach here.
For this we consider the  case of  five sinusoidal signals with frequencies falling into a `valley'  of the    noise PSD  { (Fig.\ref{Fig0}, right panel)}.
We assume the noise PSD follows an AR model (which is indeed the case here) and we  estimate $\widehat{S}_{E,AR}$ as described in Sec.\ref{sec6}.
 The question is the reliability of tests using  ${\bf \widetilde{\bf P}}\;|\;\widehat{\bf S}_{E,AR}$, and in particular how accurate is expression \eqref{PFA_param}
  { with this approach}.
  
\begin{figure}[htb!]  	\centerline{
	 \hspace{-0.2cm} \includegraphics[width= 9.5cm,height=5cm]{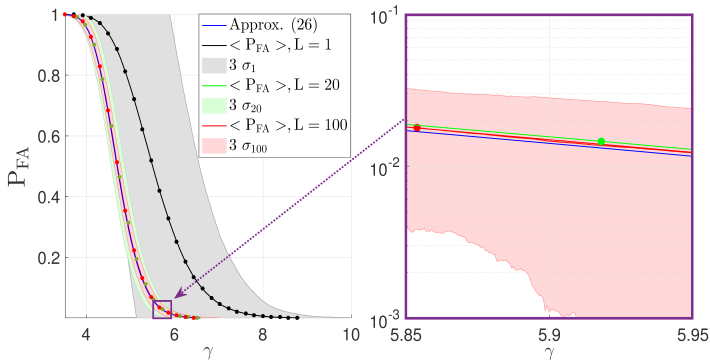}}
	\caption{ Comparison of the approximated $\rm P_{FA}$  \eqref{PFA_param} (blue curve) with true $\rm P_{FA}$ for $L=1,20$ and $100$.
	The solid lines with dots represent the average actual $\rm P_{FA}$. The shaded regions are the corresponding  $3\sigma$ enveloppes. The  right panel is a zoom in log-scale on the violet square in left panel.  { Parameters: $N_s = 5$, $\alpha_q = 0.07$ {m.s$^{-1}$}, $f_q =[5.0, 5.5, 5.75, 6, 6.5]$ mHz.}}
	\label{Fig3}
\end{figure}

 Fig.\ref{Fig3} compares, as a function of the test threshold,  the $\rm P_{FA}$ assumed by approximation \eqref{PFA_param} (blue curve) with the actual false alarm rates obtained for  $1000$ MC simulations.

 In each such simulation, an estimate $\widehat{\bf S}_{E,AR}(L)$ was obtained  with Akaike's Final Prediction Error (FPE)   \cite{Akaike_1969}  from $L$ noise time series and used to calibrate the periodogram. 
  %
%  { The choice of this criterion is not a decisive choice as, in prepapration of this study, we observed a similar behavior for all the other AR estimate criteria (CAT, RIS, AIC).}
  %
  For each such estimate, 
  %(called $\widehat{\bf S}_L$ for short hereafter), 
  the true $\rm P_{FA}$ of test 
  %  ${\rm T_C}(  {\bf \widetilde{\bf P}}\;|\; \widehat{\bf S}_{L},N_C=N_s)$
  ${\rm T_C}(  {\bf \widetilde{\bf P}}\;|\;{ \widehat{\bf S}_{E,AR}},N_C=N_s)$
   was evaluated using 100  MC simulations. The figure plots, respectively for $L=1, 20$ and $100$ the average  $\rm P_{FA}$,  respectively in black, green and red solid lines with dots. We see that \eqref{PFA_param} is accurate, \textit{on average}, only when $L$ is large.
This figure also indicates the variability of the true false alarm rate. For each value of $L$, the figure shows the $3\sigma$ dispersion of the true  $\rm P_{FA}$  w.r.t. its empirical average  (shaded regions in grey for $L=1$, green for $L=20$ and red for $L=100$). Even when $N$ is large, the true
significance levels at which such  tests are conducted can undergo wild (and in practice unknown) variations. The right panel is a zoom on the $3\sigma$ region for $L=100$.
{
For a threshold $\gamma =5.9$ for instance, the $\rm P_{FA}\approx 0.013$ from \eqref{PFA_param}. In the right panel, we see
that  the true false alarm rate for this threshold varies in reality  in the range $[0.001 \; 0.3]$.}
%We see for instance that from \eqref{PFA_param}, a threshold $\gamma =5.9$ corresponds to $\rm P_{FA}\approx 0.013$. In practice, the true false alarm rate for this threshold varies  in the range $[0.001 \; 0.3]$. 
For smaller values of $L$, the excursions of the true false alarm rates are so large that \eqref{PFA_param} is simply useless. Conclusions drawn from tests based on  parametric estimation may thus be very hazardous, even
{ when the parametric model is true and }
 when large  data sets are available for PSD estimation. 
\subsection{${\rm  T_C}({\bf \widetilde{\bf P}}|\overline{\bf P}_L,N_C)$ and ${\rm T_M}({\bf \widetilde{\bf P}}|\overline{\bf P}_L)$ vs adaptative approaches} 
 \label{SecF}
  We do not attempt to apply  ${\rm HC^\star}$ and ${\rm BJ}$ to ${\bf{P}}$ in the correlated case as this leads to the same inconsistencies as those illustrated in Sec. \ref{secC} and Fig. \ref{Fig2},   {top left panel, dotted curves.}
 We compare here the  adaptative approaches ${\rm HC^\star}$ and ${\rm BJ}$  standardized by $\overline{\bf P}_L$  as in Sec.\ref{sec_HC} with the corresponding  tests ${\rm T_M}$ and ${\rm T_C}$  for $L=1$ and $L=50$. (Note that when $N_C=1$, ${\rm T_C}$ reduces to ${\rm T_M}$, \textit{cf} \eqref{test_ch}). 
 
 { For this comparison, we consider two different cases of Keplerian signals under $\mathcal{H}_1$, which are typical of `super-Earth' planets (Fig. \ref{Fig5}):  $N_p = 1$ planet with null eccentricity in Case 1 (left column), and $N_p = 5$ eccentric planets in Case 2 (right column). 
 \begin{figure}[htb!]  \centering
	\includegraphics[width=9cm,height=7cm]{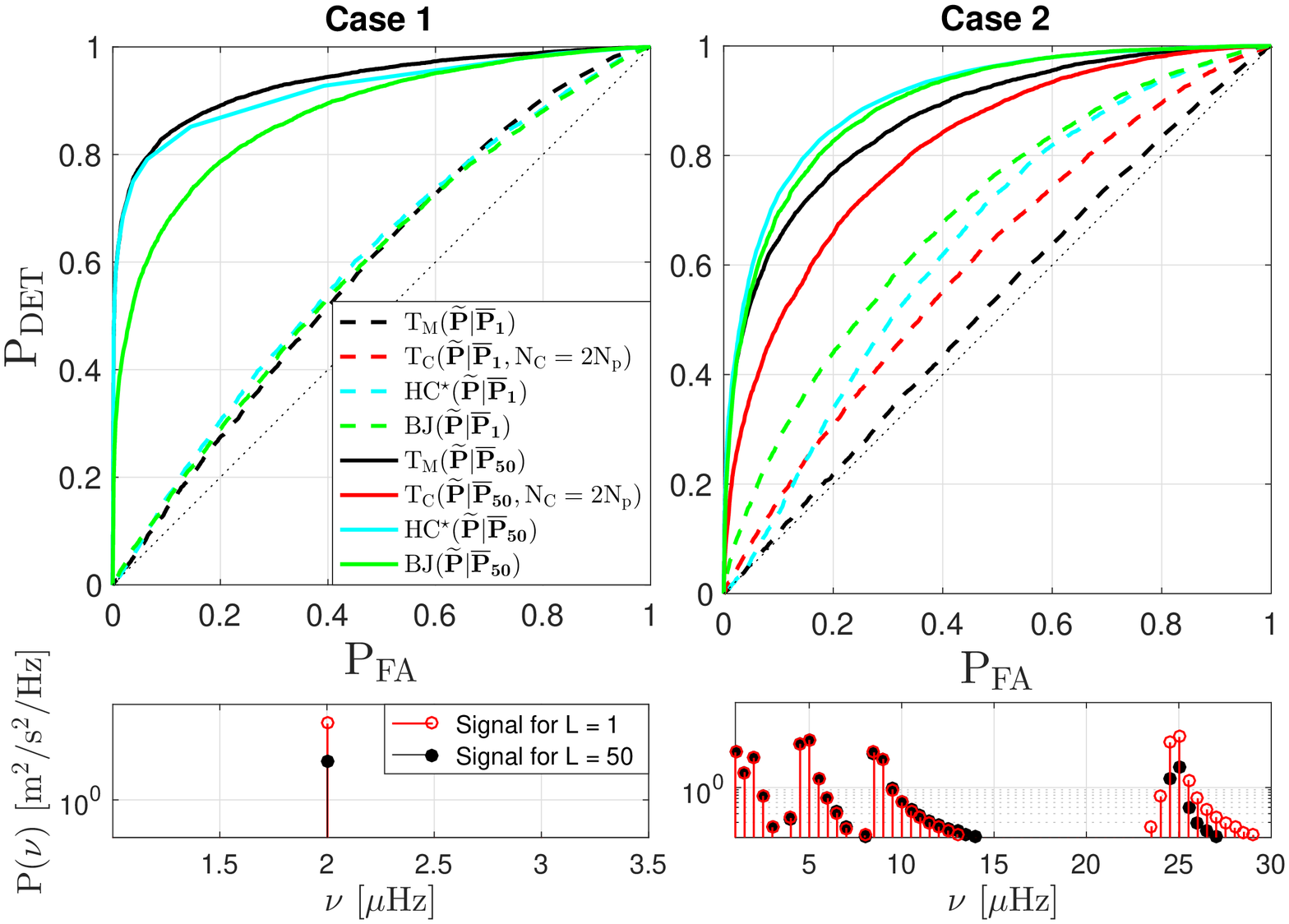}
%	\caption{ Empirical ROC curves comparing adap$N_s = [1,50]$ sines with  amplitude $\alpha_q=1 $ m.s$^{-1}$ (left) and $\alpha_q = 0.05$ m.s$^{-1}$ (right) for $L=1$ and $\alpha_q = 0.12$ m.s$^{-1}$ (left) and $\alpha_q= 0.05$  m.s$^{-1}$ (right) for $L=50$. When $N_s = 1$ the signal frequency is 5 $\rm mH_z$, while for $N_s = 50$ the signals frequencies are regularly spaced over  the whole frequency range. In  the adaptative tests,   $\alpha_0=1/2$.}
	\caption{ { Top panels: Empirical ROC curves comparing classical and adaptative tests for different signals.
	 Case 1:  $N_p = 1$, $M_p = 3.5 M_\oplus$ (for $L=1$), $T_p = 5.7813$ d (signal frequency $1/T_p$ on-grid), $e_p=0$, $\omega_p=0$, $T_0=0$, $\gamma_0= 0$. For $L = 50$, $M_p = 0.8 M_\oplus$. 
	Case 2: $N_p = 5$, $M_p =  [0.15,0.15,0.25,0.25,0.25]  M_\oplus$ (for $L=1$), $T_p = [11.21~ {\rm h}, 1.33 ~ {\rm d}, 2.45~ {\rm d}, 6.91 ~ {\rm d}, 9.25~ {\rm d}] $ (signal frequencies off-grid), $e_p=0.9$, $\omega_p=\pi$, $T_0=0$, $\gamma_0= 0$. For $L = 50$,  $M_p(T_p = 11.21~ {\rm h}) = 0.07 M_\oplus$. 
	 In  the adaptative tests,   $\alpha_0=1/2$.  Bottom panels: periodograms (logscale) of the signals under ${\mathcal{H}}_1$.}}
	\label{Fig5}
\end{figure} 

 % to compare the performance of the adaptative approaches ${\rm HC^\star}$ and ${\rm BJ}$ standardized by $\overline{\bf P}_L$  with the ${\rm T_M}$ and ${\rm T_C}$ tests .
 %Fig. \ref{Fig5} shows the results for two different numbers of planets ($N_p$) under $\mathcal{H}_1$ (case 1, left: $Np = 1$; case 2, right: $Np= 5$). 
 The lower panels illustrate the periodogram of the (noiseless) Keplerian signatures for the two cases. We see the apparition of significant harmonic{s} in the case of off-grid signal frequencies and highly eccentric  orbits. For each case, the planet masses have been adapted depending on the considered value of $L$  for a better display of the ROC curves.
 
 From $L=1$ (dashed lines) to $L=50$ (solid lines),  the performances of all  tests increase in both cases.  When the signal is extremely sparse in the Fourier domain (Case 1), ${\rm T_M}$ is more powerful than the considered adaptive approaches. This situation  changes when the spectrum is less sparse (Case 2,  compare the bottom panels), which is expected. 
 
 Results of \cite{Donoho_2004} show that for a proportion of deviations in the range $ [(\frac{N}{2}-1)^{\frac{1}{4}}\; (\frac{N}{2}-1)^{\frac{1}{2}}]$ an adaptive procedure such as HC may have better asymptotic power  than $\rm T_M$. For the case $N=1110$ considered here,  this corresponds to the range $[5\;23]$ (considering $\Omega$). The situation should be opposite  for  very sparse signals (in the range  $ [1\; 4]$ here). It turns out that the superiority of adaptive procedures is confirmed in Fig. \ref{Fig5}, 
 where the spectrum is $1$-sparse in Case 1 and about $10$-sparse in Case 2.
% for which $N_s=1$ in Case 1 and $N_s \approx 2 \times N_p = 10$ in Case 2.
Note, however, that while the theory of  \cite{Donoho_2004} may be used as a guideline for guessing the  sparsity range in which each test should work better,  we should not expect a too tight agreement with  this theory. Indeed, in the framework of \cite{Donoho_2004},
all deviations under the alternative have the same amplitude, while RV signals lead to  different amplitudes in general. Moreover, these theoretical results are asymptotical (while  $N$ is not so large here).\\
 An interesting point is the comparison of $\rm HC^\star$ and $\rm BJ$ with ${\rm T_C}({\bf \widetilde{\bf P}}|\overline{\bf P}_L,N_C)$, for which $N_C$ is a proxy for the number of significant deviations in the Fourier spectrum (here we considered $N_C=1$ in Case 1 and $N_C=10$ in Case 2). ${\rm T_C}$ represents a  kind of Oracle, which knows  in which region of the $P$-values to look at in order to `make the case' against $\mathcal{H}_0$. The right panel shows that adaptive procedures  have better power  than ${\rm T_C}$, yet without prior knowledge.}
{ \subsection{A detectability study}
\label{SecG}

A direct application of the previous results is detectability studies, which can be used for the design of observational strategies for instance. 
%A interesting case would be that of $\alpha$ Centauri Bb, but the sampling of the spectroscopic data used to claim the detection of $\alpha$ Centauri Bb \cite{Dumusque_2012} is irregular so this case cannot be addressed in the  framework of the present paper. We neverthless consider a similar detectability study below.
To illustrate this{,} we consider under  $\mathcal{H}_1$  a planet with the  parameters of $\alpha$ Centauri B's exoplanet candidate as estimated in \cite{Dumusque_2012} (see the legend of Fig.\ref{Fig6}).  As the eccentricity is supposed null, the signal can be considered very sparse in the Fourier domain. We consider the $\rm T_M$ test and a time sampling of $\Delta t = 4$ hours. It was allowed to slightly vary from one value of $N$ to another in order to guarantee
that the planet's period yields a frequency exactly on the Fourier grid, in which case the spectrum is $1$-sparse on $\Omega$ and  $\rm T_M$ is  the test that  yields the best performances.

As for the noise PSD, we used a model based on HD simulations of a star with similar spectral type as that of $\alpha$CenB. 
%There is  some mismatch between the simulated spectral type 
{ There is some mismatch between the simulated RV time series of the spectral type} (G2) and { the} true spectral type (K1) of $\alpha$CenB. Because  spectral type affects the noise properties \cite{Meunier_2016},  these results should not be considered to reflect perfectly the case of the candidate planet orbiting $\alpha$CenB. 

In Fig.\ref{Fig6}, we illustrate the feasible performance compromises $(\rm P_{DET}, \rm P_{FA})$ as a function of $N$, for three target $\rm P_{FA}$ ($0.5,\; 0.1$ and $0.01$, indicated  respectively by the dotted, solid and dashed lines) and for different sizes of available training data sets  $L$ ($ \infty, 100, 20, 5$, shown respectively in black, green, blue and red). These curves were built using the expressions   (\ref{Eq_gam}) and  (\ref{Eq20}) for the ${\rm T_M}$ test and we checked that they are accurate in separate MC simulations (not shown). In the case $L \to+\infty$, we use the fact that $ F(2,2L)  \underset{L \to+\infty}{\longrightarrow} \chi_2^2$ to calculate the theoretical $\rm P_{DET}(P_{FA})$.

The study presented in Fig.\ref{Fig6}  allows to quantify interesting facts. 
First, of course,  $\rm P_{DET}$ is larger if the allowed $\textrm{P}_{\textrm{FA}}$ is larger. Second, for a fixed   $\textrm{P}_{\textrm{FA}}$,   $\rm P_{DET}$ is larger for a larger value of $L$. Going to specific cases, we see that if a planet similar to the considered candidate was orbiting a star of the considered spectral type (G2), of which 100 training time series are available, it would require  250 days ($1500\times 4$h) of observations   with 1 point every 4 hours to guarantee a probability of detection of $0.9$ while ensuring a false alarm rate of $0.01$. This situation is indicated by the black square. With only $L=5$ training time series, the probability of detection would fall to about $0.1$, all other parameters equal (red square).
}
 \begin{figure}[htb!] \centering
	\includegraphics[width=9cm,height=5.5cm]{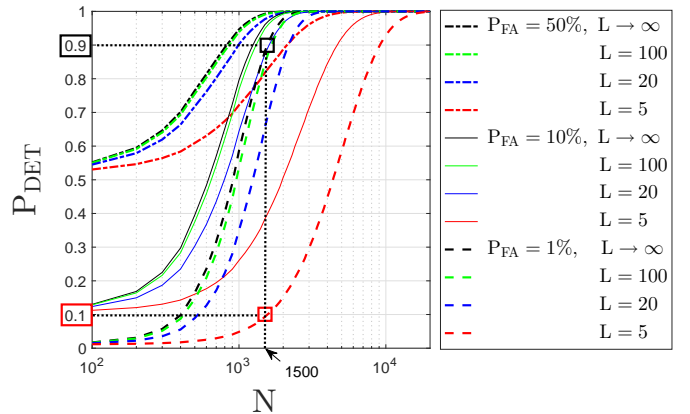}
	\caption{{Example of a detectability study for a  planet (relevant Keplerian parameters:  $K = 0.54$ m.s$^{-1}$, $T_p  = 3.23$ d and $e=0$; $M_p = 1.241 M_\oplus$, {orbit inclination of $90 \degree$}, semi-major axis $a     =   0.0425$ Astronomical Unit) orbiting a G2 type star. The plot shows the achievable $\rm P_{DET}$ for different configurations of $\textrm{P}_{\textrm{FA}}$ budgets and numbers of available training light curves.  % Fit parameters : $\hat{o} = 15$, $\widehat{\sigma^2} =2.96$; $\widehat{c_j} = [0.2076,    0.1545,    0.0329,    0.0268,   -0.0148,    0.0960,   -0.0250   -0.0048   -0.0193,   -0.0217,    0.0544,   -0.0802,   0.0832,   -0.1061,    0.0796]$ 
	}}
	\label{Fig6}
\end{figure}

\section{Summary and perspectives}
\label{conc}
This paper has  investigated the possibility of using training data sets to standardize periodograms in order to improve the control of the resulting false alarm rate. The paper first { provided}  an extended (though unavoidably selective)  overview of classical and recent techniques in sinusoid detection, with emphasis  on the problem of designing CFAR detectors for the composite hypotheses of multiple sinusoids and colored noise.
We proposed an asymptotic analysis of the periodograms statistics after standardization for a model involving an unknown number of sinusoids with unknown parameters in partially unknown colored noise. This analysis allowed
to  characterize the performances of some standardized tests in terms of false alarm and detection rates. We showed that when standardization is performed with a simple averaged periodogram as a noise PSD estimate  these tests are CFAR for all sizes of training data sets. In contrast, we pointed out that standardization based on parametric estimates of the noise PSD may present  actual false alarm rates that may be very far from the assumed ones, even for large data sets.\\
The tests we considered include classical approaches and also more recent adaptive tests designed for the rare and weak setting. For the latter tests, the standardization (by $\overline{\bf{P}}_L$) offers the same benefits as if the statistics of the noise were known {\textit{a priori}}, with the CFAR property added. 
%We also showed that such tests present power comparable to procedures for which the number of sinusoids would be known in advance.\\
{We also showed that such tests can present better power comparing to procedures for which the number of sinusoids would be known in advance.}

%   stationarity (i), regular sampling (ii) and large $N$ (iii), Sec. \ref{secIA}) 

{ In practical  situations, some of  the assumptions (i-iv) in Sec. \ref{secIA} may not be met. This can be the case  in exoplanet detection  using RV time series, owing for instance to astrophysical effects linked to magnetic activity (like spots, affecting (i)) or instrumental defects / observational constraints
 (affecting (ii) and (iii)).  Comparing theory to pratice, the present study is useful to   make feasibility studies (as in Sec. \ref{SecG}), which then describe best achievable performances (\textit{i.e.}, around ``quiet'' stars and in absence of other unmodelled perturbations) for a regular sampling. }
 %The important case of irregular sampling is discussed at the end of this Section.
   
%  We also point out that while our study and results present an obvious interest for exoplanet detection, they are sufficientl }, we mentioned other fields where they can be useful. 
 
 An important extension to the considered framework regards the case  when, for periodicity analysis, the time series is not correlated to orthogonal exponentials. This case encompasses situations where (a) the sampling is  irregular, (b)   $P$ is not evaluated on the Fourier grid (as in oversampled periodograms) (c) $P$ is modified in the form of  ``generalized periodograms'', which correlate the time series with highly redundant dictionaries of specific features 
{\cite{Bolviken_1983a,Baluev_2008,Suveges_2012, Scargle_1982,Bretthorst_2003,Thong_2004,Zechmeister_2009,Baluev_2015,Gregory_2016}}. {   In such cases,   the considered ordinates exhibit strong dependencies.
% We are currently investigating  whether asymptotic results on maxima of dependent variables \cite{Suveges_2012,Suveges_2015} can be exploited in a standardized framework.
%We are currently extending the proposed procedure based on periodogram standardization by a PSD estimate obtained through HD simulations to the case of an irregular sampling. In the case
With the additional complication of partially unknown colored noise, analytical evaluation of the false alarm rate  for the considered tests seems out of reach. %Besides,
% procedures based on an estimate of the number of ``independent frequencies'' (e.g.\cite{Horne_1986}) are unreliable in practice because
%they assume white noise.
We conducted however preliminary studies suggesting that it might be possible to obtain accurate estimates of the false alarm rate when the noise is colored,
the sampling irregular, the considered frequencies not restricted to Fourier grid and $N$ small, by
combining  the standardization proposed in the present paper  with bootstrap procedures \cite{Zoubir_1993} and maximum likelihood estimation of Generalized Extreme Values distributions' parameters \cite{Suveges_2012,Suveges_2015}.  This numerical approach, still under study, would allow to address important questions, like that of the impact of the sampling distribution on the detection performances.}
% which sampling (regular or irregular) is best from a detection viewpoint in the considered framework.}

%The use of the GEV makes  As an illustration for reviewers information, Fig \ref{fig2} shows
%in red the true false alarm rate as a function of the threshold for test $T_M$ applied to a standardized periodogram whose complex exponentials are sampled at irregular time instants, and in grey several estimates obtained by the procedure outlined above. The three bottom panels show, for three target false alarm rates,
%the true P$_{FA}$ in red (which in practice is  unknown) and the distribution of the estimates as obtained by boostrap/GEV. We see that it is still possible
%to control fairly accurately the false alarm rate in the case of irregular sampling and partially unknown colored noise. 
%
%
%
% In such case, the loss of the dependence property between the periodogram components do not allow to found its asymptotic distribution. Some empirical procedures have been derived using an estimate of the number of independent frequencies (e.g.\cite{Horne_1986}) but are unreliable in practice when the noise is colored. In this spirit, we are investigating the control of the periodogram statistics by AR estimations, MC simulations and using the Generalized Extrem Value (GEV) distribution \cite{Suveges_2012,Suveges_2015}.

% ----------------------------------------- ANNEXES  ---------------------------------------------------------
\appendices

 \section{Derivation of expressions \eqref{dist2} and \eqref{lambda}}
 \label{app1}
  
 We prove  here that for model \eqref{hyp} the periodogram is asymptotically distributed as in \eqref{dist2} with non centrality parameters as in \eqref{lambda}.
 The proof is adapted from Theorem 6.2 of \cite{Li_2014}, which considers the complex case.
  We first prove \eqref{dist2} and then turn to  \eqref{lambda}. The time series of model \eqref{hyp} can also be written as
$$
X(j)  =  \sum_{q = 1}^{N_s} \alpha_q \sin(2\pi f_q j+\varphi_q)+ E(j)= R(j)+E(j),
$$
with $R(j):=   \sum_{q=1}^{N_s}  \alpha_q \sin(2\pi  f_q j + \varphi_q)$ a deterministic part, which using Euler formulae 
can  be written as
$$
      \begin{aligned}
	R(j)  	     &=  \sum_{q=1}^{N_s} \frac{ \alpha_q}{2}  \mathrm{e}^{{\rm{i}} (\varphi_q-\frac{\pi}{2})} \mathrm{e}^{2\pi {\rm{i}} f_q j } -  \frac{ \alpha_q}{2}  \mathrm{e}^{-{\rm{i}} (\varphi_q+\frac{\pi}{2})} \mathrm{e}^{-2\pi {\rm{i}} f_q j}.
     \end{aligned}
     \vspace{-0.3cm}
 $$
$$
\vspace{-0.2cm}
\!\!\!\!\!\! \!\! \text{By introducing:}~~~~~~
	\left\{
      \begin{aligned}
      		& {\bf{f}}(\nu):=[\mathrm{e}^{{\rm{i}}2\pi \nu},\hdots, \mathrm{e}^{{\rm{i}}2\pi N\nu}]^\top,\\
		&{\bf{f}}^{+}(f_q)
		:=\mathrm{e}^{{\rm{i}} (\varphi_q-\frac{\pi}{2})}  {\bf{f}}(f_q), \\
		&{\bf{f}}^{-}(f_q):=
		\mathrm{e}^{-{\rm{i}} (\varphi_q+\frac{\pi}{2})}  {\bf{f}}(f_q),
     \end{aligned}
     \right.
$$
the time series  writes in vector form :
\begin{equation}
      \begin{aligned}
	{\bf{X}} &=   \displaystyle{\sum_{q=1}^{N_s}} \frac{\alpha_q}{2} \Big( {{{\bf{f}^{+}}(f_q)} -  {{\bf{f}^{-}}(f_q)}}  \Big)+ {\bf{E}}={\bf{R}}+ {\bf{E}} \\
  \end{aligned}
  \label{RE}
\end{equation}
and its discrete Fourier transform (DFT) $y_k$ at frequency $\nu_k$ { is}
$$
	y_k=\frac{1}{N}  {\bf{f}}^H(\nu_k){\bf{X}} = \frac{1}{N}  {\bf{f}}^H(\nu_k){\bf{R}} + \frac{1}{N}  {\bf{f}}^H(\nu_k){\bf{E}}.
$$
The DFT is composed of a deterministic part, ${\mu_k}$, and a stochastic part, ${\epsilon}_k$, defined as
\begin{equation}
	{\mu_k}:= \frac{1}{N} {\bf{f}}^H(\nu_k){\bf{R}}\quad{\textrm{and}}\quad {\epsilon}_k:= \frac{1}{N} {\bf{f}}^H(\nu_k)\bf{E}.
	\label{mu}
\end{equation}
 Because $E$ is a zero mean Gaussian process, the random variable $y_k=\mu_k+\epsilon_k$ is Gaussian with mean $\mu_k$ and variance $\sigma^2_k$. This is a complex variable for all Fourier frequencies except for $\nu_0$ and $\nu_{\frac{N}{2}}$ because $ {\bf{f}}(\nu_0)$ and $ {\bf{f}}(\nu_{\frac{N}{2}})$ are real.

The distribution of the periodogram requires to investigate the variance of $\epsilon_k$ which, with \eqref{mu}, writes:
$$
      \begin{aligned}
      &N \textrm{var}{\;\epsilon_k}=\frac{1}{N}\mathbb{E}\; \left( {\bf{f}}^H(\nu_k){\bf{E}}{\bf{E}}^\top  {\bf{f}}(\nu_k)\right) \\
	&=\frac{1}{N}\displaystyle{\sum_{t,s=1}^N}r_E(t-s)\mathrm{e}^{-{\rm{i}}2\pi \nu_k(t-s)} \\
	&= \frac{1}{N}\displaystyle{\sum_{|u|<N}}(N-|u|)r_E(u)\mathrm{e}^{-{\rm{i}}2\pi \nu_ku}\\
	&= S_{E}(\nu_k)-\!\!\!\displaystyle{\sum_{|u|<N}}\frac{|u|}{N}r_E(u)\mathrm{e}^{-{\rm{i}}2\pi \nu_ku}
 -\!\!\!\displaystyle{\sum_{|u|\geq N}}r_E(u)\mathrm{e}^{-{\rm{i}}2\pi \nu_ku}\\
	 & =S_{E}(\nu_k) + {\mathcal{O}}(r_N),
   \end{aligned}
$$
since, using { the absolutely summable autocorrelation function $\sum_u |r_E(u)| < \infty$} and the dominated convergence theorem we have
\vspace{-0.4cm}
\begin{equation}
	r_N:=\sum_u \min(1,\frac{|u|}{N}) |r_E(u)| \to 0 \textrm{\;\;as\;\;} N \to \infty.
	\vspace{-0.2cm}
	\label{lerho} 
\end{equation}
Hence, for all Fourier frequencies,
 \begin{equation}
	 \sigma^2_k:= \textrm{var}\;\epsilon_k= N^{-1} S_{E}(\nu_k) + {\mathcal{O}}(N^{-1} r_N).
	 \label{lavar}
 \end{equation}
 
 By Lemma 12.2.1(b) of \cite{Li_2014}, we obtain 
$$
  |y_k|^2/\sigma^2_k \sim   \left\{
      \begin{aligned}
		&{  \frac{1}{2}} \chi^2_{2,2\frac{|\mu_k|^2}{\sigma^2_k}},  ~~~~ \forall k ~ \in ~ \Omega,\\
		&  \chi^2_{1,\frac{|\mu_k|^2}{\sigma^2_k}}, ~ \text{ for } ~ k = 0, \frac{N}{2}.
	 \end{aligned}
	  \right.
$$

Hence,  for the periodogram this implies
\begin{equation} 
\begin{aligned} 
P(\nu_k | H_1)& = N|y_k|^2 = N\sigma^2_k(|y_k|^2/\sigma^2_k)\\
& \sim
  \left\{         
      \begin{aligned}
	 &\frac{1}{2}\rho^{-1}_k S_E(\nu_k) \chi_{2, 2\rho_k \gamma_k}^2 ,~ \forall k ~ \in ~ \Omega,  \\
	 &\rho^{-1}_k S_E(\nu_k)  \chi^2_{1,\rho_k \gamma_k }, ~ \text{ for } ~ k = 0, \frac{N}{2},
      \end{aligned}
    \right.
      \end{aligned}
      \vspace{-0.2cm}
       \label{leP}
\end{equation}
where 
\vspace{-0.2cm}
\begin{equation}
\rho_k:=S_E(\nu_k)/(N\sigma^2_k) \quad{\textrm{and}}\quad\gamma_k:=N|\mu_k|^2/S_E(\nu_k).
\label{gamma}
\end{equation} 
With \eqref{lavar}, we see that
\begin{equation}
	\rho_k=\frac{S_E(\nu_k)}{S_E(\nu_k)+{\mathcal{O}(r_N)}}=1+{\mathcal{O}(r_N)}
	 \label{lerho2}
\end{equation}
 for all Fourier frequencies.
 Owing to \eqref{lerho}, an approximated distribution of  \eqref{leP} can be obtained by neglecting
  the $\mathcal{O}(r_N)$ in { \eqref{lerho2}}. The distribution \eqref{dist2} follows by noting
\begin{equation}
	\hspace{-3mm}
	\lambda_k:=2\gamma_k\;\; \text{for}\;\; k ~ \in ~ \Omega\quad \text{and}\quad \lambda_k:=\gamma_k \;\; \text{for}  \;\; k = 0, \frac{N}{2}.
\label{fact}
\end{equation} 
 We now turn to the computation of the non centrality parameters. We have from \eqref{RE}, \eqref{mu} and \eqref{gamma}
\begin{equation} \small 
	  \begin{aligned}
	\gamma_k &= \frac{N}{S_E(\nu_k)} |  \frac{1}{N} {\bf{f}}^H(\nu_k)\bf{R}|^2\\
	 &= \frac{1}{NS_E(\nu_k)}  \Big| \sum_{j =1}^{N}   \sum_{q=1}^{N_s}  \frac{\alpha_q}{2} \Big( \mathrm{e}^{{\rm{i}} (\varphi_q-\frac{\pi}{2})} \mathrm{e}^{2\pi {\rm{i}} (f_q-\nu_k) j } \! \\
	 &~~~~~~~~~~~~~~~~~~~~~~~~~~\hdots  - \!\mathrm{e}^{-{\rm{i}} (\varphi_q+\frac{\pi}{2})} \mathrm{e}^{-2\pi {\rm{i}} (f_q+\nu_k) j} \Big) \Big|^2\!. \\
	 \end{aligned}
 	\label{gamma3}
\end{equation}
{ Introducing the Dirichlet Kernel (\textit{cf} Lemma 12.1.3 of \cite{Li_2014}):\\ }
$$
	D_N(\nu) := \frac{1}{N}  \sum_{j=1}^{N}  \mathrm{e}^{{\rm{i}} 2\pi \nu j} = \frac{\sin(N\pi \nu)}{N\sin(\pi \nu)} \mathrm{e}^{{\rm{i}} (N+1)\pi \nu},
     \label{eq_Dn}
$$
	and the corresponding  Fej\'er kernel (or spectral window)
\begin{equation}
	K_N(\nu) := |D_N(\nu)|^2 =  \Bigg( \frac{\sin(N\pi \nu)}{N\sin(\pi \nu)} \Bigg)^2,
     \label{eq_Kn}
\end{equation}
expression \eqref{gamma3} becomes
{
$$ 
 \small
  \begin{aligned}
 \gamma_k 
	%&=  \frac{1}{N S_E(\nu_k)} \Big| N \sum_{q=1}^{N_s}  \frac{\alpha_q}{2} \Big( D_N(f_q -\nu_k)  \mathrm{e}^{{\rm{i}} (\varphi_q-\frac{\pi}{2})} \!-\! D_N(f_q +\nu_k) \mathrm{e}^{-{\rm{i}} (\varphi_q+\frac{\pi}{2})}\! \Big)\! \Big|^2 \\
				&=  \frac{N}{4 S_E(\nu_k)} \Big| \sum_{q=1}^{N_s}  \alpha_q \Big( D_N(f_q -\nu_k)  \mathrm{e}^{{\rm{i}} (\varphi_q-\frac{\pi}{2})} \\
				&~~~~~~~~~~~~~~~~~~~~~~~~~~~~~~~~~~~~\hdots  - D_N(f_q +\nu_k) \mathrm{e}^{-{\rm{i}} (\varphi_q+\frac{\pi}{2})}\! \Big)\! \Big|^2. \\
 \end{aligned}
 $$
 }
 This equation can also be written as :
\begin{equation} 
	 \gamma_k =  \frac{N}{4 S_E(\nu_k)} \Big| \sum_{q=1}^{N_s}  \alpha_q z_q(\nu_k) \Big|^2,
	 \label{gamma2}
 \end{equation} 
 with
 \begin{equation} 
  \begin{aligned}
	z_q(\nu_k) &:= D_N(f_q -\nu_k)  \mathrm{e}^{{\rm{i}} (\varphi_q-\frac{\pi}{2})} \!-\! D_N(f_q +\nu_k) \mathrm{e}^{-{\rm{i}}(\varphi_q+\frac{\pi}{2})} \\
	        &=  x_{+}  \mathrm{e}^{{\rm{i}} \theta_+ } -  x_{-}  \mathrm{e}^{ {\rm{i}} \theta_{-} }\\
		&=  ( x_{+}  \cos\theta_+ -   x_{-}  \cos\theta_{-} )+ {\rm{i}}  ( x_{+} \sin\theta_+ -   x_{-}\sin\theta_{-} ),
 \end{aligned}
 \label{zp}
\end{equation}
where
$$
\left\{
  \begin{aligned}
 	x_{+}&=  x_{+}\left({\nu_k},q\right)  :=  \frac{\sin(N\pi (f_q - \nu_k))}{N\sin(\pi (f_q - \nu_k))},  \\
 	x_{-}&=  x_{-}\left({\nu_k},q\right)    :=  \frac{\sin(N\pi (f_q + \nu_k))}{N\sin(\pi (f_q + \nu_k))} , \\
	\theta_+ &= \theta_+\left({\nu_k},q\right) :=  + [ (N+1)\pi (f_q - \nu_k)+ (\varphi_q-\frac{\pi}{2})], \\
	\theta_- &= \theta_-\left({\nu_k},q\right)  :=  -[ (N+1)\pi (f_q + \nu_k)+ (\varphi_q+\frac{\pi}{2})]. \\	
 \end{aligned}
 \right.
$$
The modulus $\kappa_q$ of $z_q$ may be written as  
  \begin{equation} 
 	\begin{aligned}
		 & \kappa_q  :=  | z_q | = \left( x_{+}^2 +  x_{-}^2 - 2 x_{+}  x_{-}  \cos{(\theta_+-\theta_{-}) }\right)^{\frac{1}{2}},
 	\end{aligned}
	\vspace{-0.4cm}
 	 \label{eq_keppa}	
 \end{equation} 

$$
\text{with}~~~~~~~~~
\left\{
 \begin{aligned}
 	x_{+}^2 &= K_N(f_q - \nu_k), \\
	x_{-}^2 &= K_N(f_q+ \nu_k) , \\
	\theta_+-\theta_{-} &=   2 \pi (N+1) f_q + 2 \varphi_q,
 \end{aligned}
 \right.
$$
\begin{equation} 
\text{and}~~~~~~~~~~~~~~~~~ \theta_q   := \angle\; z_q ,\; \theta_q\; \in\;]-\pi, \pi],
  \label{eq_theta}
  \end{equation} 
  the phase of $z_p$ obtained from the real and imaginary parts of \eqref{zp} \cite{Kasana_2005}.
With these notations, it  is easy to show that 
 $$
 \Big| \sum_{q=1}^{N_s} z_q \Big|^2 = \sum_{q=1}^{N_s} \Big[ \kappa_q^2 + 2  \kappa_q \sum_{\ell = q+1}^{N_s} \kappa_\ell \cos(\theta_q-\theta_\ell)  \Big] ~~ \text{for $N_s>1$}.
 $$
Consequently, the  expression of the $\{\gamma_k\}$ in \eqref{gamma2} becomes 
\begin{equation} 
	 \begin{aligned}
		&\gamma_k \!=\!  \frac{N}{4 S_E(\nu_k)} \! \! \sum_{q=1}^{N_s} \!\! \Big[  \alpha_q^2 \kappa_q ^2 \!\!+\! 2\alpha_q \kappa_q\!\! \sum_{\ell = q+1}^{N_s}\!\! \! \alpha_\ell \kappa_\ell  \!\cos(\theta_q-\theta_\ell )\Big]	\\
	\end{aligned}
 \label{eq_50}	
\end{equation} 
 and  the non centrality parameters $\{\lambda_k\}$ of \eqref{lambda} follow from \eqref{fact},with $\kappa_q$ and $\theta_q$  given by \eqref{eq_keppa} and \eqref{eq_theta}.  $\QEDB$

Note that if all signal frequencies $\{f_p\}$ fall on the Fourier frequency grid, the crossed term in \eqref{eq_50} vanish owing to the orthogonality
of the Fej\'er kernels centered at different signal frequencies. In this case, expression \eqref{eq_50} precisely reduces to expression given in Remark 6.6 of \cite{Li_2014}.

We finally wish to mention that the expression of the non centrality parameters is erroneously reported in exp. (5) of  \cite{Sulis_2016a} (sign error and crossed terms missing).

% ------------------------------------------------------------------------------------------------------------------------------------------------
%  PDET FOR TC TEST
% ------------------------------------------------------------------------------------------------------------------------------------------------

 \section{Derivation of expression  (\ref{lagrosse}) }
   \label{app6} 
   Let $K$ denote the number of  ordinates of  ${\bf \widetilde{\bf P}}\;|\;\overline{\bf P}_L$  larger than $\gamma$  under $\mathcal{H}_1$,
 and $p_i := \textrm{Pr}\; ( K=i\;|\;{\cal{H}}_1)$.  From the definition \eqref{test_ch}, we have:
  \begin{equation}  
\begin{aligned} 
	&\!\!\rm{P_{DET}}({T}_{C}( {\bf \widetilde{\bf P}}|\overline{\bf P}_L),\gamma,N_C) \!:=\! \textrm{Pr}\; ( {T_C  ({\bf \widetilde{\bf P}}|\overline{\bf P}_L,N_C)}\!\! > \!\gamma | {\cal{H}}_1\!)\! \\
	&=  \textrm{Pr}\; ( K\geq N_C\;|\;{\cal{H}}_1) \\
	&= 1- { \sum_{i=0}^{N_C-1} } \textrm{Pr}\; ( K=i\;|\;{\cal{H}}_1) =  1- \displaystyle{ \sum_{i = 0}^{N_C-1}} p_i.
\end{aligned}
\label{pfab}
 \end{equation}
 
\noindent Owing to \eqref{dist_H1} each  ordinate $({\bf \widetilde{\bf P}}|\overline{\bf P}_L)_i:=\frac{\widetilde{P}(\nu_i)}{\overline{ P}_L(\nu_i)}$ has probability $1-\Phi_{F_{\lambda_i}}(\gamma)$ to be larger than $\gamma$. These variates can be considered approximately independent but not i.i.d. Hence, the variable
$K$ is not binomially distributed (as it is under ${\cal{H}}_0$) and the probabilities $\{p_i\}$ require further investigation.  
We proceed by induction.  In the following, all probabilities are under $\mathcal{H}_1$. The first probability  $p_0$ can simply be approximated as 
$$
 \begin{aligned} 
		&p_0=  \textrm{Pr} \left\{\bigcap_{ k=1}^{\frac{N}{2}-1}
	   ({\bf \widetilde{\bf P}}\;|\;\overline{\bf P}_L)_k \leq \gamma   \right\}  \approx \displaystyle{ \prod_{k =1}^{ \frac{N}{2}-1}} \Phi_{F_{\lambda_{k}}}. \\
\end{aligned} 
$$
The probability $p_1=\text{Pr}(K=1)$ is similarly
$$
 \begin{aligned} 
 	p_1 &= \textrm{Pr} {\displaystyle{\bigcup_{k=1}^{\frac{N}{2}-1} }}\left\{  ({\bf \widetilde{\bf P}}\;|\;\overline{\bf P}_L)_k> \gamma      \bigcap_{ j\neq k}
	   ({\bf \widetilde{\bf P}}\;|\;\overline{\bf P}_L)_j \leq \gamma   \right\}   \\
	&\approx \displaystyle{ \sum_{k = 1}^{\frac{N}{2}-1} \Big[  ( 1 -  \Phi_{F_{\lambda_{k}}} ) \displaystyle{ \prod_{\substack{j =1,\\ j \neq k}}^{ \frac{N}{2}-1}} \Phi_{F_{\lambda_{j}}}  \Big] }. \\
\end{aligned} 
$$	
To generalize further, denote by ${\Omega}^{(i)}$ one particular combination of $i$ indices taken in $\Omega$ and $\overline{\Omega}^{(i)}:=\Omega \backslash \Omega^{(i)}$ the set of remaining indices. Let $\{\Omega^{(i)}_1,\hdots , \Omega^{(i)}_i\}$ (resp. $\{\overline{\Omega}^{(i)}_1,\hdots, \overline{\Omega}^{(i)}_{\frac{N}{2}-1-i}\}$) denote the indices in  two such combinations, and let   $\Omega^i$  be the set of all the  $\{{\Omega}^{(i)}\}$. 
With these notations we obtain  for $i>1$ :
\begin{equation}
 \begin{aligned} 
 	p_i &= \textrm{Pr} \!\!\! {\displaystyle{\bigcup_{ {\Omega}^{(i)} \in {\Omega}^{i} } }}
	\left\{   \bigcap_{ k=1}^i	   ({\bf \widetilde{\bf P}}\;|\;\overline{\bf P}_L)_{\Omega^{(i)}_k } > \gamma
	   \bigcap_{ k'=1}^{\frac{N}{2}-1-i}	   ({\bf \widetilde{\bf P}}\;|\;\overline{\bf P}_L)_{\overline{\Omega}^{(i)}_{k'}}  \leq \gamma \right\}   \\	
	&\approx  \displaystyle{  \sum_{\Omega^{(i)} } \prod_{k=1}^{ i  } } \Big(\!1\! - \Phi_{F_{\lambda_{\Omega^{(i)}_k}}\!\!^{\!\!\!\!\!\!\!\!\!(\gamma, 2,2L)}} \!\!\Big) \displaystyle{  \prod_{k'=1}^{ \frac{N}{2}-1-i }}\Phi_{F_{\lambda_{ \overline{\Omega}^{(i)}\!_{\!\!\!\!\!\!k'} }}{(\gamma, 2,2L)} }.
\end{aligned} 
\label{pfac}
\end{equation}	
 Expression  (\ref{lagrosse}) follows by combining \eqref{pfab} and \eqref{pfac}.  $\QEDB$
	
% ----------------------------------------- ACKNOWLEDGMENT  -----------------------------------------
\section*{Acknowledgement}
The GOLF instrument onboard SoHO is a cooperative effort of scientists, engineers, and technicians, to whom we are indebted. SoHO is a project of international collaboration between ESA and NASA.
% ----------------------------------------- BIBLIO  -----------------------------------------
\bibliographystyle{IEEEbib}
\bibliography{maBiblio,Biblio_IEEE} 

\begin{thebibliography}{10}

\bibitem{Schuster_1898}
A.~{Schuster},
\newblock ``{On the investigation of hidden periodicities},''
\newblock {\em J. Geophys. Res.}, vol. 3, pp. 13, 1898.

\bibitem{Brillinger_1981}
D.R. {Brillinger},
\newblock {\em Time Series : Data Analysis and Theory},
\newblock Holden Day, San Francisco, 1981.

\bibitem{Bartlett_1950}
M.S. {Bartlett},
\newblock ``Periodogram analysis and continuous spectra,''
\newblock {\em Biometrika}, vol. 37, pp. 1--16, 1950.

\bibitem{Grenander_1957}
U.~{Grenander} and M.~{Rosenblatt},
\newblock {\em Statistical analysis of stationary time series},
\newblock John Wiley and Sons, 1957.

\bibitem{Priestley_1981}
M.B. {Priestley},
\newblock {\em Spectral Analysis and Time Series},
\newblock Academic Press, San Diego, 1981.

\bibitem{Brockwell_1991}
P.J. Brockwell and R.A. Davis,
\newblock {\em Time series : theory and methods},
\newblock Springer, 1991.

\bibitem{Bloomfield_2000}
P.~{Bloomfield},
\newblock {\em Fourier Analysis of Time Series},
\newblock Wiley-Intersci., 2000.

\bibitem{Stoica_2005}
P.~{Stoica} and R.~{Moses},
\newblock {\em Spectral analysis of signals},
\newblock Prentice Hall, 2005.

\bibitem{Quinn_2001}
B.~G. {Quinn} and E.J. {Hannan},
\newblock {\em The Estimation and Tracking of Frequency},
\newblock Cambridge Univ., 2001.

\bibitem{Pepe_2014b}
F.~{Pepe} et~al.,
\newblock ``{Instrumentation for the detection and characterization of
  exoplanets},''
\newblock {\em Nature}, vol. 513, pp. 358--366, 2014.

\bibitem{Batalha_2014}
N.M. {Batalha},
\newblock ``{Exploring exoplanet populations with NASA's Kepler Mission},''
\newblock {\em Proc. Nat. Acad. Sci.}, vol. 111, pp. 12647--12654, Sept. 2014.

\bibitem{Auvergne_2009}
M.~{Auvergne} et~al.,
\newblock ``{The CoRoT satellite in flight: description and performance},''
\newblock {\em A\&A}, vol. 506, pp. 411--424, 2009.

\bibitem{Rauer_2014}
H.~{Rauer} et~al.,
\newblock ``{The PLATO 2.0 mission},''
\newblock {\em Experimental Astronomy}, vol. 38, pp. 249--330, 2014.

\bibitem{Fischer_2014}
D.A. {Fischer} et~al.,
\newblock ``{Exoplanet Detection Techniques},''
\newblock {\em Protostars and Planets VI}, pp. 715--737, 2014.

\bibitem{Perryman_2011}
M.~{Perryman},
\newblock {\em The exoplanet handbook},
\newblock Cambridge Univ., 2011.

\bibitem{Dumusque_2012}
X.~{Dumusque} et~al.,
\newblock ``{An Earth-mass planet orbiting {$\alpha$} Centauri B},''
\newblock {\em Nature}, vol. 491, pp. 207--211, 2012.

\bibitem{Hatzes_2013}
A.~{Hatzes},
\newblock ``{The Radial Velocity Detection of Earth-mass Planets in the
  Presence of Activity Noise: The Case of {$\alpha$} Centauri Bb},''
\newblock {\em ApJ}, vol. 770, pp. 133, 2013.

\bibitem{Rajpaul_2016}
V.~{Rajpaul} et~al.,
\newblock ``{Ghost in the time series: no planet for Alpha Cen B},''
\newblock {\em MNRAS}, vol. 456, pp. L6--L10, 2016.

\bibitem{Bigot_2011}
L.~{Bigot} et~al.,
\newblock ``{The diameter of the CoRoT target HD 49933},''
\newblock {\em A\&A}, vol. 534, no. 3, 2011.

\bibitem{Fisher_1929}
R.A. {Fisher},
\newblock ``{Tests of Significance in Harmonic Analysis},''
\newblock {\em Proc. R. Soc. London, Ser. A}, vol. 125, pp. 54--59, 1929.

\bibitem{Chiu_1989}
S.T. Chiu,
\newblock ``Detecting periodic components in a white gaussian time series,''
\newblock {\em J. R. Stat. Soc. Series B}, vol. 51, no. 2, pp. 249--259, 1989.

\bibitem{Shimshoni_1971}
M.~{Shimshoni},
\newblock ``On fisher's test of significance in harmonic analysis,''
\newblock {\em Geophys. J. R. Astronom. Soc.}, pp. 373--377, 1971.

\bibitem{Kay_1999}
S.~Kay,
\newblock ``Adaptive detection for unknown noise power spectral densities,''
\newblock {\em IEEE Trans. Signal Process}, vol. 47, no. 1, pp. 10--21, 1999.

\bibitem{Li_2014}
T.H. {Li},
\newblock {\em Time series with mixed spectra},
\newblock CRC Press, 2014.

\bibitem{Whittle_1952}
P.~Whittle,
\newblock ``The simultaneous estimation of a time series harmonic components
  and covariance structure,''
\newblock {\em Trabajos de Estadistica}, vol. 3, no. 1-2, pp. 43--57, 1952.

\bibitem{Bartlett_1955}
M.S. Bartlett,
\newblock ``An introduction to stochastic processes,''
\newblock {\em Quart. J. R. Meteorological Soc.}, vol. 81, no. 350, pp. 650,
  1955.

\bibitem{Siegel_1980}
A.~Siegel,
\newblock ``Testing for periodicity in a time series,''
\newblock {\em J. Amer. Stat. Assoc.}, vol. 75, no. 370, pp. 345--348, 1980.

\bibitem{Bolviken_1983a}
E.~B{\"o}lviken,
\newblock ``New tests of significance in periodogram analysis,''
\newblock {\em Scandinavian J. Stat.}, vol. 10, no. 1, pp. 1--9, 1983.

\bibitem{Bolviken_1983b}
E.~B{\"o}lviken,
\newblock ``The distribution of certain rational functions of order statistics
  from exponential distributions,''
\newblock {\em Scandinavian J. Stat.}, vol. 10, no. 2, pp. 117--123, 1983.

\bibitem{Sachs_1993}
R.~Von~Sachs,
\newblock ``Estimating the spectrum of a stochastic process in the presence of
  a contaminating signal,''
\newblock {\em IEEE Trans. Signal Process.}, vol. 41, no. 1, pp. 323, 1993.

\bibitem{Sachs_1994}
R.~Von~Sachs,
\newblock ``Peak-insensitive non-parametric spectrum estimation,''
\newblock {\em J. Time Series Anal.}, vol. 15, no. 4, pp. 429--452, 1994.

\bibitem{Bhansali_1979}
R.J. Bhansali,
\newblock ``A mixed spectrum analysis of the lynx data,''
\newblock {\em J. R. Stat. Soc. Ser. A}, , no. 142, pp. 199--209, 1979.

\bibitem{Truong_1990}
B.~Truong-Van,
\newblock ``A new approach to frequency analysis with amplified harmonics,''
\newblock {\em J. R. Stat. Soc. Ser. B}, vol. 52, pp. 203--221, 1990.

\bibitem{Quinn_1991}
B.G. {Quinn} and J.M. {Fernandes},
\newblock ``A fast efficient technique for the estimation of frequency,''
\newblock {\em Biometrika}, vol. 78, no. 3, pp. 489--497, 1991.

\bibitem{Quinn_1999}
B.G. {Quinn},
\newblock ``A fast efficient technique for the estimation of frequency:
  Interpretation and generalisation,''
\newblock {\em Biometrika}, vol. 86, pp. 213--220, 1999.

\bibitem{Kavalieris_1994}
L.~{Kavalieris} and {Hannan} E.J.,
\newblock ``Determining the number of terms in a trigonometric regression.,''
\newblock {\em J. Time Series Analysis}, vol. 15, pp. 613--625, 1994.

\bibitem{Hannan_1961}
E.J. Hannan,
\newblock ``Testing for a jump in the spectral function,''
\newblock {\em J. R. Stat. Soc. Ser. B}, vol. 23, no. 2, pp. 394--404, 1961.

\bibitem{Nicholls_1967}
D.F. Nicholls,
\newblock ``Estimation of the spectral density function when testing for a jump
  in the spectrum,''
\newblock {\em Austr. J. Stat.}, vol. 9, pp. 103--108, 1967.

\bibitem{Chiu_1990}
S.T. {Chiu},
\newblock ``{Peak-insensitive parametric spectrum estimation},''
\newblock {\em Stochastic Processes and their Applications}, vol. 35, pp.
  121--140, 1990.

\bibitem{Gryca_1998}
J.K. Gryca,
\newblock ``Detection of multiple sinusoids buried in noise via balanced model
  truncation,''
\newblock {\em IEEE Instrum. Meas. Conf.}, pp. 1353--1358, 1998.

\bibitem{White_1999}
L.B. White,
\newblock ``Detection of sinusoids in unknown coloured noise using ratios of ar
  spectrum estimates,''
\newblock {\em Proc. Inform., Decision and Control}, pp. 257--262, 1999.

\bibitem{Lu_2005}
N.~{Lu} and D.~{Zimmerman},
\newblock ``Testing for directional symmetry in spatial dependence using the
  periodogram,''
\newblock {\em J. Stat. Planning and Inference}, vol. 129, pp. 369--385, 2005.

\bibitem{Liavas_1998}
A.P. {Liavas} et~al.,
\newblock ``A periodogram-based method for the detection of steady-state
  visually evoked potentials.,''
\newblock {\em IEEE Trans. Biom. Eng.}, vol. 45, no. 2, pp. 242--248, 1998.

\bibitem{Zheng_2012}
C.~Zheng,
\newblock ``Detection of multiple sinusoids in unknown colored noise using
  truncated cepstrum thresholding and local signal-to-noise-ratio,''
\newblock {\em Applied Acoust.}, pp. 809--816, 2012.

\bibitem{Nadler_2011}
B.~{Nadler} and A.~{Kontorovich},
\newblock ``Model selection for sinusoids in noise: Statistical analysis and a
  new penalty term,''
\newblock {\em IEEE Trans. Signal Process.}, vol. 59, no. 4, pp. 1333--1345,
  2011.

\bibitem{Koen_2015a}
C.~{Koen},
\newblock ``The analysis of indexed astronomical time series - xi. the
  statistics of oversampled white noise periodograms,''
\newblock {\em MNRAS}, vol. 449, no. 1, pp. 1098--1105, 2015.

\bibitem{Koen_2015b}
C.~{Koen},
\newblock ``The analysis of indexed astronomical time series - xii. the
  statistics of oversampled fourier spectra of noise plus a single sinusoid,''
\newblock {\em MNRAS}, vol. 453, no. 2, pp. 1793--1798, 2015.

\bibitem{Tuomi_2012}
M.~{Tuomi} et~al.,
\newblock ``{Signals embedded in the radial velocity noise. Periodic variations
  in the {$\tau$} Ceti velocities},''
\newblock {\em A\&A}, vol. 551, pp. A79, 2012.

\bibitem{Donoho_2004}
D.~{Donoho} and J.~{Jin},
\newblock ``{Higher criticism for detecting sparse heterogeneous mixtures},''
\newblock {\em Ann. Stat.}, 2004.

\bibitem{Ingster_2010}
Y.I. {Ingster} et~al.,
\newblock ``{Detection boundary in sparse regression},''
\newblock {\em Electron. J. Statist.}, vol. 4, pp. 1476--1526, 2010.

\bibitem{Walther_2011}
G.~{Walther},
\newblock ``{The Average Likelihood Ratio for Large-scale Multiple Testing and
  Detecting Sparse Mixtures},''
\newblock {\em IMS Collections : From Probability to Stat. and Back:
  High-Dimensional Models and Processes}, vol. 9, pp. 317--326, 2012.

\bibitem{Moscovich_2016}
A.~{Moscovich} et~al.,
\newblock ``{On the exact Berk-Jones statistics and their p-value
  calculation},''
\newblock {\em Electron. J. Stat.}, vol. 10, pp. 2329--2354, 2016.

\bibitem{Gontscharuk_2016}
V.~{Gontscharuk} et~al.,
\newblock ``Goodness of fit tests in terms of local levels with special
  emphasis on higher criticism tests,''
\newblock {\em Bernoulli}, vol. 22, no. 3, pp. 1331--1363, 2016,.

\bibitem{Hall_2010}
P.~{Hall} and J.~{Jin},
\newblock ``Innovated higher criticism for detecting sparse signals in
  correlated noise,''
\newblock {\em Ann. Stat.}, vol. 38, pp. 1686--1732, 2010.

\bibitem{Sulis_2016a}
S.~{Sulis}, D.~{Mary}, and L.~{Bigot},
\newblock ``Using hydrodynamical simulations of stellar atmospheres for
  periodogram standardization: Application to exoplanet detection,''
\newblock in {\em IEEE ICASSP}, March 2016, pp. 4428--4432.

\bibitem{Gupta_2001}
S.K. {Gupta} et~al.,
\newblock ``{UPSO three channel fast photometer},''
\newblock {\em Bulletin Astron. Soc. of India}, vol. 29, pp. 479--486, 2001.

\bibitem{Sulis_2016c}
S.~{Sulis}, D.~{Mary}, and L.~{Bigot},
\newblock ``{Overcoming the stellar noise barrier for the detection of telluric
  exoplanets: an approach based on hydrodynamical simulations},''
\newblock in {\em EAS Publications Series}, D.~{Mary} et~al., Eds., Sept. 2016,
  vol. 78-79, pp. 247--274.

\bibitem{Proakis_1996}
J.G. {Proakis} and D.G. {Manolakis},
\newblock {\em Digital Signal Processing},
\newblock Prentice-Hall, 1996.

\bibitem{Abramowitz_1972}
M.~{Abramowitz} et~al.,
\newblock {\em Spectral Analysis and Time Series},
\newblock Dover Publications, 1972.

\bibitem{David_2003}
H.A. {David} and H.N. {Nagaraja},
\newblock {\em Order Statistics, 3rd Ed.},
\newblock Wiley, 2003.

\bibitem{Kay_1998}
S.M. {Kay},
\newblock {\em Fundamentals of Statistical signal processing. Vol II :
  Detection theory.},
\newblock Prentice-Hall, Inc, 1998.

\bibitem{Quinn_1986}
B.~G. {Quinn},
\newblock ``Testing for the presence of sinusoidal components,''
\newblock {\em J. Appl. Probability}, vol. 23, pp. 201--210, 1986.

\bibitem{Schwarzenberg_1998}
A.~{Schwarzenberg-Czerny},
\newblock ``{The distribution of empirical periodograms: Lomb--Scargle and PDM
  spectra},''
\newblock {\em MNRAS}, vol. 301, pp. 831--840, 1998.

\bibitem{Aittokallio_2001}
T.~{Aittokallio} et~al.,
\newblock ``{Testing for Periodicity in Signals: An Application to Detect
  Partial Upper Airway Obstruction during Sleep,},''
\newblock {\em J. Theoretical Medicine}, vol. 3, no. 4, pp. 231--245, 2001.

\bibitem{Guitierrez_2009}
J.~{Guti{\'e}rrez-Soto} et~al.,
\newblock ``{Low-amplitude variations detected by CoRoT in the B8IIIe star HD
  175869},''
\newblock {\em A\&A}, vol. 506, pp. 133--141, 2009.

\bibitem{Aldor_2013}
S.~{Aldor-Noiman} et~al.,
\newblock ``The power to see: A new graphical test of normality.,''
\newblock {\em Am. Stat.}, vol. 68, no. 4, pp. 318--318, 2013.

\bibitem{Mary_2014}
D.~Mary and A.~Ferrari,
\newblock ``A non-asymptotic standardization of binomial counts in higher
  criticism,''
\newblock in {\em Inform. Theory (ISIT), IEEE Int. Symp.}, June 2014, pp.
  561--565.

\bibitem{Kaplan_2014}
D.M. {Kaplan} and M.~{Goldman},
\newblock ``True equality (of pointwise sensitivity) at last: a dirichlet
  alternative to {K}olmogorov-{S}mirnov inference on distributions.,''
\newblock {\em Tech. report}, 2014.

\bibitem{Gontscharuk_2014}
V.~{Gontscharuk} et~al.,
\newblock ``The intermediates take it all: Asymptotics of higher criticism
  statistics and a powerful alternative based on equal local levels.,''
\newblock {\em Biom. J.}, vol. 57, no. 1, pp. 159--180, 2014.

\bibitem{Li_2015}
J.~{Li} and D.~{Siegmund},
\newblock ``Higher criticism: $p$-values and criticism,''
\newblock {\em Ann. Statist.}, vol. 43, no. 3, pp. 1323--1350, 2015.

\bibitem{Berk_1979}
R.H. {Berk} and D.H. {Jones},
\newblock ``Goodness-of-fit test statistics that dominate the {K}olmogorov
  statistics,''
\newblock {\em Z. Wahrscheinlichkeit.}, vol. 47, pp. 47--59, 1979.

\bibitem{Moscovich_2015b}
A.~{Moscovich-Eiger} and B.~{Nadler},
\newblock ``{Fast calculation of boundary crossing probabilities for Poisson
  processes},''
\newblock {\em ArXiv e-prints (V.3)}, 2015.

\bibitem{Fan_2003}
J.~{Fan} and Q.~{Yao},
\newblock {\em Nonlinear Time Series-Nonparametric and Parametric Methods},
\newblock Springer-Verlag New York, 2003.

\bibitem{Turkman_1984}
K.~F. {Turkman} and A.~M. {Walker},
\newblock ``On the asymptotic distributions of maxima of trigonometric
  polynomials with random coefficients,''
\newblock {\em Advances in Applied Probability}, vol. 16, no. 4, pp. 819--842,
  1984.

\bibitem{Davis_1999}
R.A {Davis} and T.~{Mikosch},
\newblock ``The maximum of the periodogram of a non-gaussian sequence,''
\newblock {\em Ann. of Prob.}, vol. 27, pp. 522--536, 1999.

\bibitem{Rife_1974}
D.~{Rife} and R.~{Boorstyn},
\newblock ``Single tone parameter estimation from discrete-time observations,''
\newblock {\em IEEE Trans. Inf. Theory}, vol. 20, no. 5, pp. 591--598, 1974.

\bibitem{Quinn_1994}
B.G. {Quinn} and P.~J. {Kootsookos},
\newblock ``Threshold behavior of the maximum likelihood estimator of
  frequency,''
\newblock {\em IEEE Trans. Signal Process.}, vol. 42, no. 11, pp. 3291--3294,
  1994.

\bibitem{Akaike_1969}
H.~Akaike,
\newblock ``Fitting autoregressive models for prediction,''
\newblock {\em Ann. Inst. Stat. Math.}, vol. 21, no. 1, pp. 243--247, 1969.

\bibitem{Akaike_1974}
H.~Akaike,
\newblock ``{A new look at the statistical model identification},''
\newblock {\em IEEE Trans. Automatic Control}, vol. 19, no. 6, pp. 716--723,
  1974.

\bibitem{Parzen_1975}
E.~Parzen,
\newblock ``{Multiple time series: determining the order of approximating
  autoregressive schemes},''
\newblock {\em Tech. Report}, , no. 23, pp. 716--723, 1975.

\bibitem{Hannan_1979}
E.J. {Hannan} and B.G. {Quinn},
\newblock ``The determination of the order of an autoregression,''
\newblock {\em J. R. Stat. Soc. Ser. B}, vol. 41, pp. 190--195, 1979.

\bibitem{Rissanen_1984}
J.~Rissanen,
\newblock ``Universal coding, information prediction and estimation,''
\newblock {\em IEEE Trans. Inf. Theory}, , no. 30, pp. 629--636, 1984.

\bibitem{Boardman_2002}
A.~{Boardman} et~al.,
\newblock ``A study on the optimum order of autoregressive models for heart
  rate variability,''
\newblock {\em Physio. Meas.}, vol. 23, pp. 325, 2002.

\bibitem{Garcia_2005}
R.A. {Garcia} et~al.,
\newblock ``{Global solar Doppler velocity determination with the GOLF/SoHO
  instrument},''
\newblock {\em A\&A}, vol. 442, pp. 385--395, 2005.

\bibitem{Meunier_2016}
N.~{Meunier} et~al.,
\newblock ``{Variability of stellar granulation and convective blueshift with
  spectral type and magnetic activity. I. K and G main sequence stars},''
\newblock {\em ArXiv e-prints}, 2016.

\bibitem{Baluev_2008}
R.V. {Baluev},
\newblock ``{Assessing the statistical significance of periodogram peaks},''
\newblock {\em MNRAS}, vol. 385, pp. 1279--1285, 2008.

\bibitem{Suveges_2012}
M.~{S{\"u}veges},
\newblock ``Extreme-value modelling for the significance assessment of
  periodogram peaks,''
\newblock {\em MNRAS}, vol. 440, no. 3, pp. 2099--2114, 2014.

\bibitem{Scargle_1982}
J.D. {Scargle},
\newblock ``{Studies in astronomical time series analysis. II - Statistical
  aspects of spectral analysis of unevenly spaced data},''
\newblock {\em ApJ}, vol. 263, pp. 835--853, 1982.

\bibitem{Bretthorst_2003}
G.L. {Bretthorst},
\newblock ``{Frequency Estimation And Generalized Lomb-Scargle
  Periodograms.},''
\newblock {\em Stat. Challenges in Astronomy}, pp. 309--329, 2003.

\bibitem{Thong_2004}
T.~{Thong} et~al.,
\newblock ``Lomb-wech periodogram for non-uniform sampling,''
\newblock {\em Proc. 26th Annu. Int. Conf. IEEE EMBS}, 2004.

\bibitem{Zechmeister_2009}
M.~{Zechmeister} and M.~{K{\"u}rster},
\newblock ``{The generalised Lomb-Scargle periodogram. A new formalism for the
  floating-mean and Keplerian periodograms},''
\newblock {\em A\&A}, vol. 496, pp. 577--584, 2009.

\bibitem{Baluev_2015}
R.V. {Baluev},
\newblock ``{Keplerian periodogram for Doppler exoplanet detection: optimized
  computation and analytic significance thresholds},''
\newblock {\em MNRAS}, vol. 446, pp. 1478--1492, 2015.

\bibitem{Gregory_2016}
P.C. Gregory,
\newblock ``An apodized kepler periodogram for separating planetary and stellar
  activity signals,''
\newblock {\em MNRAS}, vol. 458, pp. 2604--2633, 2016.

\bibitem{Zoubir_1993}
A.M. {Zoubir},
\newblock ``{Bootstrap: theory and applications},''
\newblock in {\em SPIE Conf.}, F.T. {Luk}, Ed., 1993, vol. 2027, pp. 216--235.

\bibitem{Suveges_2015}
M.~S{\"u}veges et~al.,
\newblock ``A comparative study of four significance measures for periodicity
  detection in astronomical surveys,''
\newblock {\em MNRAS}, vol. 450, no. 2, pp. 2052--2066, 2015.

\bibitem{Kasana_2005}
H.S. {Kasana},
\newblock {\em Complex variables : theory and applications. 2nd Edition},
\newblock Prentice-Hall of India, 2005.

\end{thebibliography}

\end{document}